\documentclass[11pt]{amsart}
\usepackage{hyperref}
\usepackage{tikz}
\usepackage[active]{srcltx}
\usepackage{xypic}
\usepackage{graphicx}

\usepackage[active]{srcltx}
\usepackage{epsfig,amsmath}
\usepackage[english]{babel}
\usepackage{amsmath,amsfonts,amssymb,amscd,color,epsfig,amsthm}

\newtheorem{newthm}{Theorem}
\newtheorem{theorem}{Theorem}[section]
\newtheorem{lemma}[theorem]{Lemma}
\newtheorem{proposition}[theorem]{Proposition}
\newtheorem{corollary}[theorem]{Corollary}

\newtheorem{definition}[theorem]{Definition}

\theoremstyle{remark}
\newtheorem{example}[theorem]{\bf Example}

\theoremstyle{plain}

\numberwithin{equation}{section}

\newcommand{\wt}{\widetilde}

\def\AAA{{\cal A}}

\def\FFF{{\cal F}}

\def\JJJ{{\cal J}}

\def\PPP{{\cal P}}
\def\QQQ{{\cal Q}}
\def\RRR{{\cal R}}

\def\WWW{{\cal W}}

\def\g{\gamma}
\def\G{\Gamma}

\def\R{\mbox{$\mathbb R$}}

\def\C{\mbox{$\mathbb C$}}
\def\T{\mbox{$\mathbb T$}}

\def\D{\mathbb D}

\def\Z{\mbox{$\mathbb Z$}}

\def\lv{ \left(\begin{matrix} }
 \def\rv{\end{matrix}\right)}

\def\cal{\mathcal}

\def\dw{{\dw}}

\newcommand{\mylabel}[1]{\label{#1}}

\newcommand{\REFEQN}[1] { \begin{equation}\mylabel{#1} }
\newcommand{\ENDEQN}{\end{equation}}
\newcommand{\REFTHM}[1] { \begin{theorem}\mylabel{#1} }
\newcommand{\ENDTHM}{\end{theorem}}
\newcommand{\REFNTH}[1] { \begin{newthm}\mylabel{#1} }
\newcommand{\ENDNTH}{\end{newthm}}
\newcommand{\REFPROP}[1]{\begin{proposition}\mylabel{#1} }
\newcommand{\ENDPROP}{\end{proposition} }
\newcommand{\REFLEM}[1]{\begin{lemma}\mylabel{#1} }
\newcommand{\ENDLEM}{\end{lemma} }
\newcommand{\REFCOR}[1]{\begin{corollary}\mylabel{#1} }
\newcommand{\ENDCOR}{\end{corollary} }

\def\pf{postcritically finite }

\def\ov{\overline}
\def\pf{postcritically finite }

\def\T{{\mathbb T}}

\usepackage{tikz}
\usetikzlibrary{arrows}
\tikzstyle{every picture}=[> = to]

\tikzset{cdlabel/.style={execute at begin node=$\scriptstyle,execute at end node=$}}
\tikzset{implication/.style={double equal sign distance, -implies}}
\tikzset{biimplication/.style={double equal sign distance, implies-implies}}

\title{The core entropy for polynomials of higher degree}
\author{Yan Gao,  Giulio Tiozzo}

\thanks{Mathematical School, Sichuan University, Chengdu, P. R. China. \url{email:gyan@scu.edu.cn}}

\thanks{Department of Mathematics, University of Toronto, Canada. \url{email:tiozzo@math.utoronto.ca}}

\date{}

\begin{document}
\maketitle

\begin{abstract}
As defined by W. Thurston, the \emph{core entropy} of a polynomial is the entropy
of the restriction to its Hubbard tree. For each $d \geq 2$, we study the core entropy as a
function on the parameter space of polynomials of degree $d$, and prove it varies continuously
both as a function of the combinatorial data and of the coefficients of the polynomials.
This confirms a conjecture of W. Thurston.
\end{abstract}

\section{Introduction}

A classical way to measure the topological complexity of a dynamical system is its \emph{entropy}.
In particular, to each real polynomial map $f$ one can associate the topological entropy of $f$ as
a dynamical system on the real line \cite{MT}.

If $f : \mathbb{C} \to \mathbb{C}$ is a complex polynomial map, then the real line is no longer invariant, and it becomes
less obvious to define a notion of entropy for $f$.
However, in the case $f$ is \emph{postcritically finite} (i.e., the forward orbits of the critical points are finite) then there is a canonical tree
inside the complex plane, known as the \emph{Hubbard tree} $H_f$, which is invariant under forward iteration \cite{DH}.

\begin{figure} \label{F:global}
\includegraphics[width = 0.95 \textwidth]{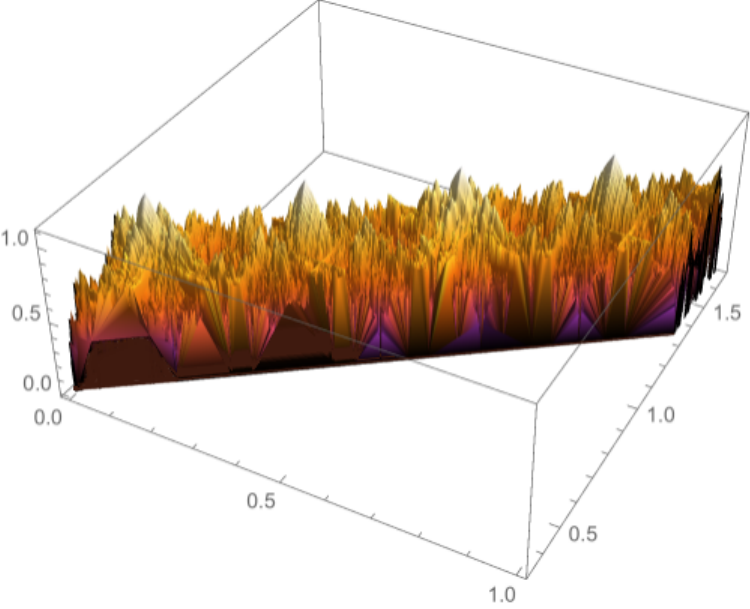}
\caption{The core entropy of cubic polynomials as a function of the primitive majors.}
\end{figure}

In order to generalize the theory of entropy to complex polynomials, W. Thurston defined the \emph{core entropy} of $f$
as the topological entropy of the restriction of $f$ to its Hubbard tree:
$$h(f) := h_{top}(f \mid_{H_f})$$
Thurston conjectured that the core entropy is a continuous function of the polynomial.
For quadratic polynomials, this was proven by \cite{Ti} and \cite{DS}.

In this paper, we generalize this result by developing the theory of core entropy for polynomials of any degree $d \geq 2$, and proving that
it varies continuously over parameter space.

\begin{figure}
\includegraphics[width = 0.9 \textwidth]{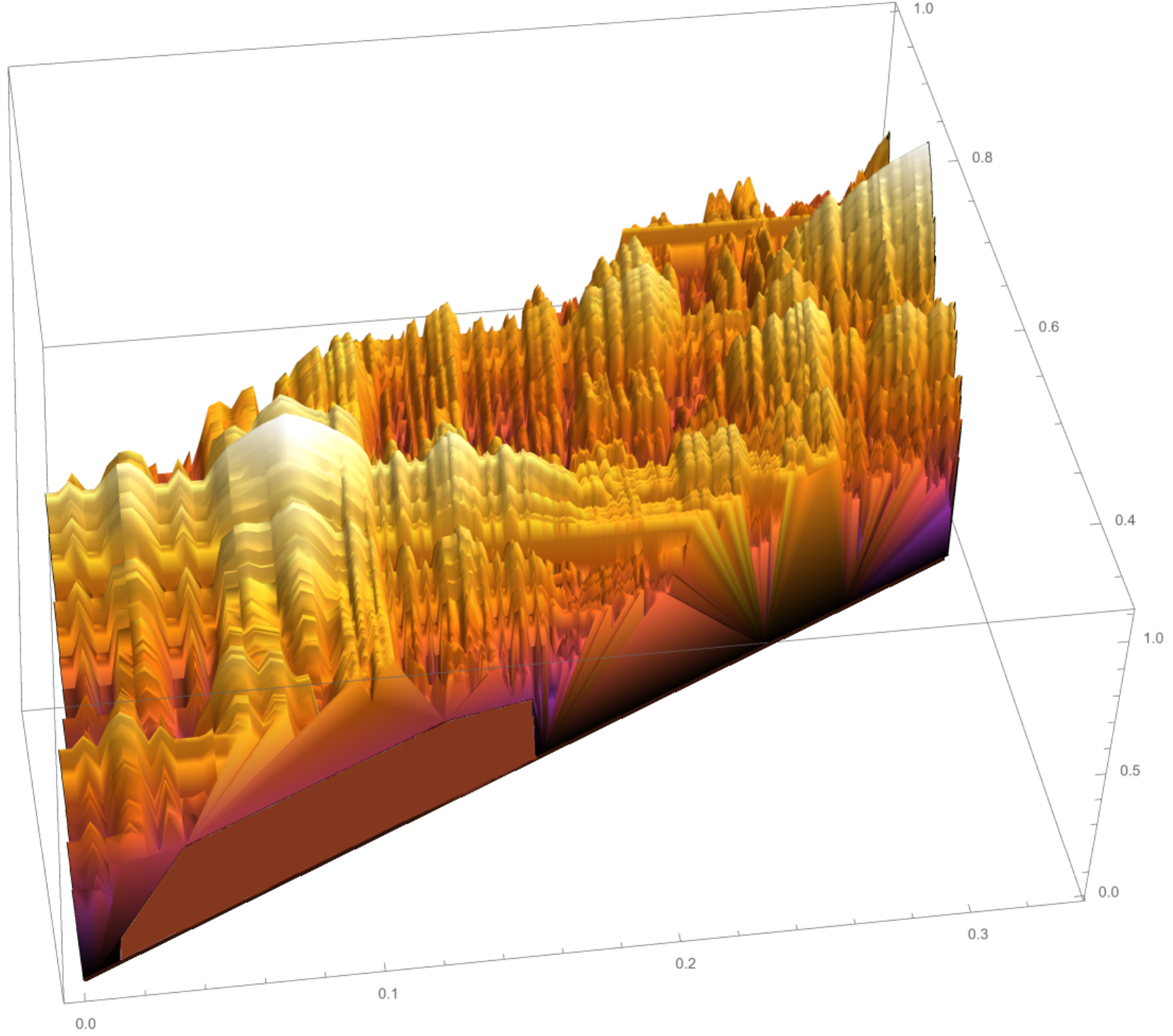}
\caption{The core entropy of cubic polynomials on the space of primitive majors. Because of the symmetries of parameter space,
it is enough to restrict attention to the domain $0 \leq a \leq 1/3$, $1/3 \leq b \leq 1$.}
\end{figure}

In order to describe the global topology of the space of polynomials of a given degree $d$, W. Thurston defined the space PM$(d)$ of \emph{primitive majors} of degree $d$ \cite{TG}.
The combinatorial parameter space PM$(d)$ generalizes to higher degree the circle at infinity for the Mandelbrot set: it is always compact and has interesting topology  (for instance, PM$(d)$ is a $K(B_d, 1)$ where $B_d$ is the braid group \cite{TG}; see also Section \ref{S:prim}).
Rational primitive majors are associated to \pf polynomials: essentially, one records the major leaf for each critical point.
Thus, we can assign to each primitive major $m$ the \emph{core entropy} $h(m)$ of the associated polynomial.
We prove that the core entropy extends to a continuous function on the combinatorial parameter space:

\begin{theorem}\label{theorem:main}
The core entropy function $h(m)$  extends to a continuous function on the set {\rm PM}$(d)$ of primitive majors of degree $d$.
\end{theorem}

The previous result is purely combinatorial, as the core entropy can be computed from the combinatorial data of the rational primitive major $m$.
In the second part of the paper, we will address the dependence of the core entropy as a function of the coefficients of the polynomial.

Let us fix $d \geq 2$, and let us consider the space $\mathcal{P}_d$ of monic, centered, polynomials of degree $d$. A polynomial is \emph{centered}
if the barycenter of its roots is the origin. Inside $\mathcal{P}_d$ sits the \emph{connectedness locus} $\mathcal{M}_d$,
i.e. the set of polynomials with connected Julia set.  While it is conjectured that the Mandelbrot set is locally connected \cite{DH}, it is known that the connectedness locus for cubics is not locally connected (\cite{La}, \cite{Mil2}, \cite{EY}).

We say that a sequence $(f_n)$ of polynomials in $\mathcal{P}_d$ converges to a polynomial $f$ if the coefficients of $f_n$
converge to the coefficients of $f$. We will show that the dependence is continuous not only
as a function of the combinatorial parameters but also as a function of the coefficients:

\begin{theorem}\label{theorem:main2}
Let $d \geq 2$. Then the core entropy $h(f)$ is a continuous function on the space of \pf polynomials of degree $d$.
\end{theorem}

To put our result in context, let us note that in going from the quadratic to the general degree case, the parameter space has complex dimension $> 1$, 
and there are several cases in dynamics where continuity fails in higher dimension: for instance, 
the action of the mapping class group on Teichm\"uller space $\mathcal{T}(S)$ extends continuously to the boundary if and only if $\textup{dim}_\mathbb{C} \mathcal{T}(S) = 1$ \cite{KT}. Similarly, Thurston's pullback map for polynomials fails to extend continuously to the Thurston boundary in higher dimension (\cite{Se}, \cite{BEKP}).

\subsection{History}
The study of topological entropy for real, quadratic polynomials goes back to the seminal work of Milnor-Thurston \cite{MT},
who proved that it depends continuously and monotonically on the parameter. Alternative proofs are also given in \cite{Do}, \cite{Ts}.

Two types of generalization of these results are possible: on the one hand, for complex quadratic polynomials generalizations of the real line
are \emph{veins} in the Mandelbrot set $\mathcal{M}$. In fact, it is known that the core entropy is monotone, increasing from the center of the Mandelbrot set
to the tips (\cite{Li}, \cite{Ti1}, \cite{Ze}). Thus, the core entropy function is intimately related to the topological structure of the Mandelbrot set:
in fact, sublevels of the entropy function can be used to define wakes in $\mathcal{M}$.

On the other hand, a considerable amount of work has gone into understanding entropy for real polynomials of higher degree.
In particular, it was conjectured by Milnor that the entropy for real polynomials of a given degree $d$ is also monotone, in the sense that isentropic curves in parameter space are connected. This was proven by Milnor-Tresser for cubics \cite{MiTr} and by Bruin-van Strien for general degrees \cite{BvS}.

The present paper is one of the first attempts to study the core entropy for higher-degree polynomials which are not real.
Note that some of the techniques used in the quadratic case do not generalize, as we cannot use the vein structure of the Mandelbrot set.
In fact, our proof of continuity does not rely at all on the understanding of the topology of the connectedness locus: on the other hand,
we believe the core entropy may be a useful tool to define and investigate the hierarchical structure of the connectedness locus,
which is much less understood than in the quadratic case.

Finally, the core entropy is related by the simple formula (see \cite{Ti1}, \cite{BS})
$$\textup{H.dim }B(f) = \frac{h(f)}{\log d}$$
to the dimension of the set $B(f)$ of biaccessible angles, which has been
studied extensively, e.g. in \cite{Za}, \cite{Zd}, \cite{Sm}, \cite{BS}, \cite{MS}.

\subsection{The techniques}

In our approach to core entropy, we avoid using the geometry or topology of parameter space, which can get extremely complicated, 
and rather we describe the combinatorics of polynomials through the use of laminations as in \cite{ThLam}, \cite{TG}. 

In the degree $2$ case, there is only one critical point, hence the combinatorics of a polynomial is captured by 
one \emph{minor leaf}. In the higher degree case, however, there are more than one critical points, hence to each polynomial one associates a finite lamination, called a \emph{critical portrait}, and to that a graph and a growth rate. 

The main difficulty to overcome with respect to the quadratic case is that the space of critical portraits is partitioned
into \emph{strata} given by the multiplicity and mutual position of the critical points, and each stratum is not closed (see also Section \ref{S:strata}). 
A recursive formula for the number of strata for each degree has been obtained by Tomasini \cite{To}. 
For instance, a sequence of polynomials with distinct critical points can converge to a polynomial with critical points of higher multiplicity, and the number 
of critical leaves used to represent it can change. To overcome this, we develop the theory of \emph{weak critical markings}, 
and we prove that \emph{all} possible limit critical portraits have the same entropy, even though they belong to different strata.

In more detail, we use an algorithm devised by Thurston \cite{TG} in order to compute the core entropy without the need to understand the topology of the Hubbard tree. In fact, Gao Yan \cite{G} proved that the algorithm yields the correct value of the core entropy for all postcritically finite polynomials of any degree (see also \cite{Ju}).

\begin{enumerate}
\item
We define the \emph{growth rate} for each critical portrait $\xi$ as follows. Given a primitive major, we construct an infinite graph $\Gamma_\xi$ (called a \emph{wedge}) whose vertices are the pairs of postcritical angles, and whose edges are given by the action of the dynamics on the space of arcs between postcritical points.

\item
We then associate to this infinite graph its \emph{growth rate} $r(\xi)$ by considering the growth rate of the number of closed paths in the graph:
$$r(\xi) := \limsup_{n \to \infty} \left( \#\{\textup{closed paths in }\Gamma_\xi \textup{ of length }n \} \right)^{1/n}$$
and prove that such number is the zero of a convergent power series $P(t)$, which we call \emph{spectral determinant}.

\item
We prove that the growth rate $r(\xi)$ depends continuously on the space of critical portraits.

\item
For rational critical portraits $\xi$, we prove that the growth rate is related to the core entropy $h(\xi)$ given by Thurston's algorithm (Lemma \ref{lem:entropy-continuity}), namely
$$h(\xi) = \log r(\xi).$$
\end{enumerate}
This establishes Theorem \ref{theorem:main} and concludes the combinatorial part of the paper.

To get the second main result (Theorem \ref{theorem:main2}), let us consider a sequence $f_n \to f$ of \pf polynomials of degree $d$.
We study the variation in the landing points of the external rays corresponding to the portrait.

\begin{enumerate}
\item Following Poirier \cite{Poi}, we construct for each \pf polynomial a rational critical portrait, known as \emph{critical marking}.
Let us pick for each $f_n$ a critical marking, which we denote as $\Theta_n$.

\item
We then study the possible limits of the sequence $(\Theta_n)$: it turns out that limits of critical markings for $f_n$ need not be critical markings for $f$, since more rays than expected can land on the same critical point (see Example \ref{E:not-conv}). For this reason, we introduce the more general notion of \emph{weak critical marking}, and prove that (Proposition \ref{Julia-type-converge}):
 $$\textup{each limit }\Theta_\infty \textup{ of the sequence }(\Theta_n)\textup{ is a weak critical marking of }f.$$
This is the analytic part of the paper, as it requires controlling the convergence of landing rays as parameters change. In fact, the parts of the marking
associated to Fatou critical points and to Julia critical points need to be dealt with separately.

\item We use \cite{G} to conclude that Thurston's algorithm also gives the correct value of core entropy for weak critical markings, hence $h(\Theta_n) = h(f_n)$ and also $h(\Theta_\infty) = h(f)$.

\item By the first combinatorial part, $h(\Theta_n) \to h(\Theta_\infty)$, hence
$$h(f_n) = h(\Theta_n) \to h(\Theta_\infty) = h(f)$$
which completes the proof of Theorem \ref{theorem:main2}.
\end{enumerate}

\subsection{The space of primitive majors} \label{S:prim}

Before delving into the proofs, we explore the structure and topology of the set PM$(d)$ of primitive majors  in degree $2$ and $3$.
See also \cite{BH}, and \cite{TG} for the general case.

\begin{figure}
\fbox{\includegraphics[width = 0.9 \textwidth]{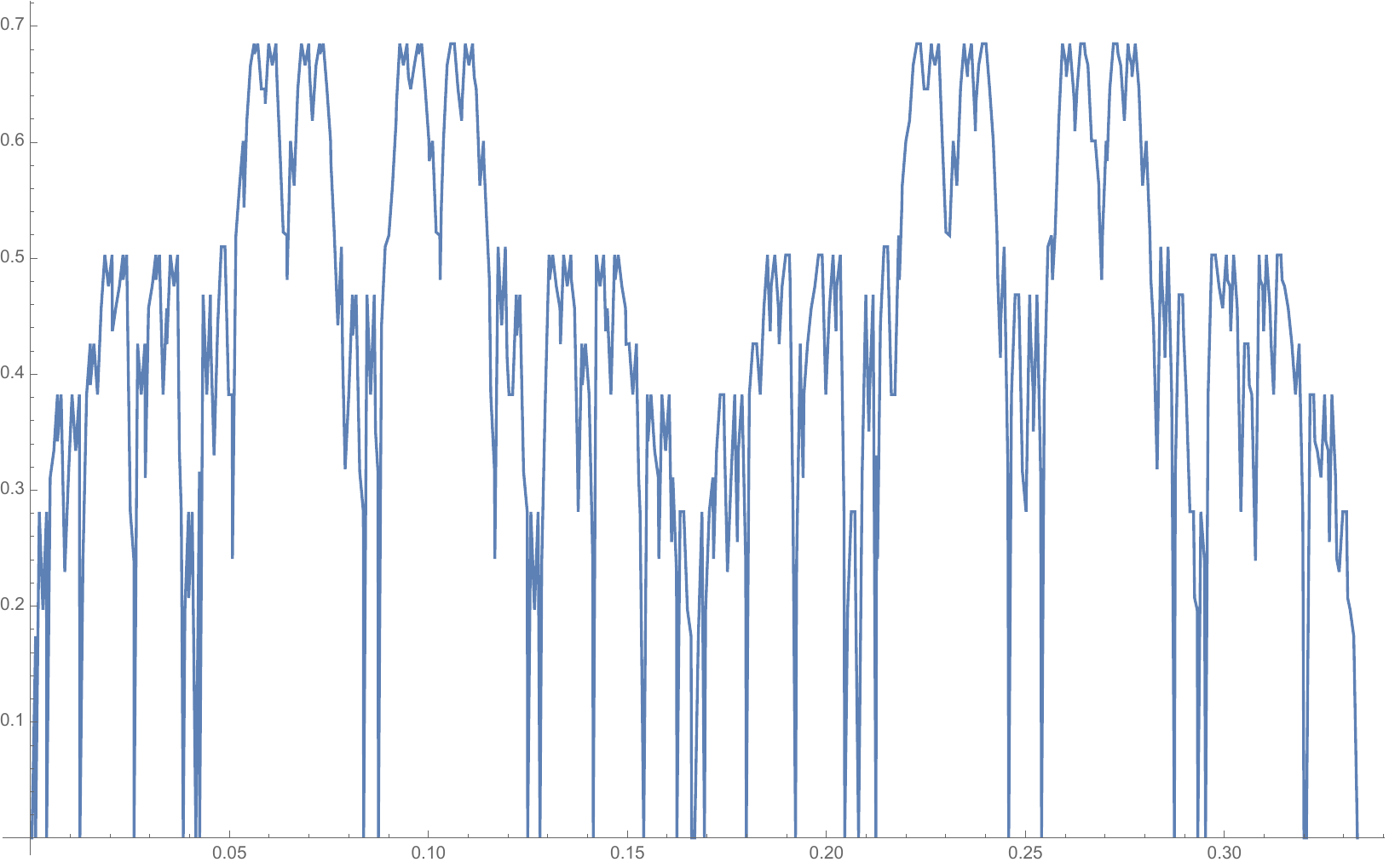}}
\caption{The core entropy for \emph{unicritical} cubic polynomials. The critical portraits are of the form $\{(a, a + 1/3, a + 2/3)\}$ with $0 \leq a \leq 1/3$. The maxima reach the value $\log 2$.}
\label{F:uni}
\end{figure}

\begin{figure}
\fbox{\includegraphics[width = 0.9 \textwidth]{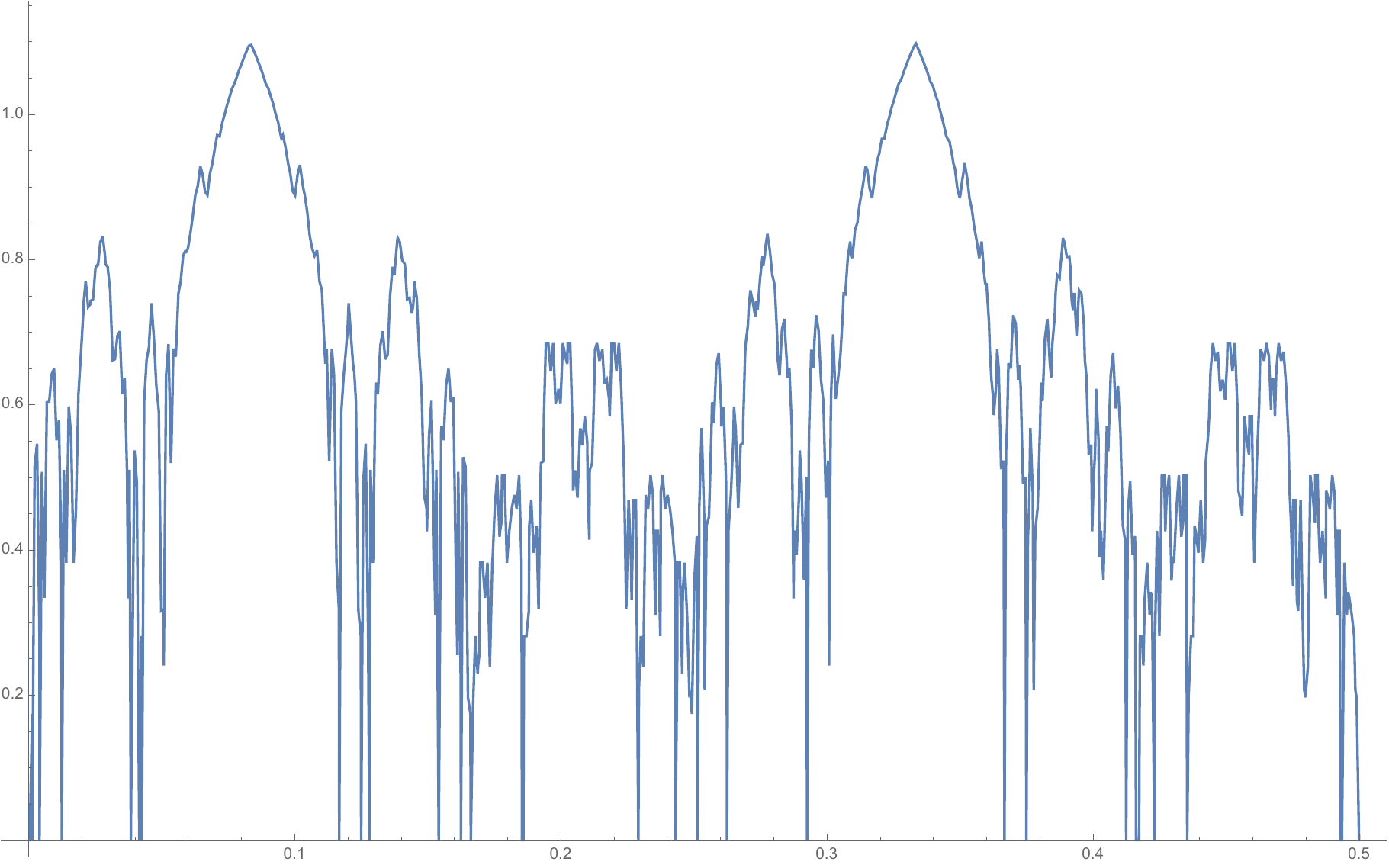}}
\caption{The core entropy for \emph{symmetric} cubic polynomials. The critical portraits are of the form $\{ (a, a + 1/3), (a + 1/2, a + 5/6)\}$ with $0 \leq a \leq 1/2$. There are two maxima for $a = 1/12$ and $a = 1/3$. These are also global maxima over the whole space PM$(3)$, and the entropy equals $\log 3$. }
 \label{F:sym}
\end{figure}

A \emph{critical portrait} of degree $d$ is defined as a collection
\[m=\{\ell_1,\ldots,\ell_s\}\]
of leaves and ideal polygons in $\ov{\D}$ fullfilling the following conditions:
\begin{enumerate}
\item any two distinct elements $\ell_k$ and $\ell_l$ either are disjoint or intersect at one point on $\partial\D$;
\item the vertices of each $\ell_k$ are identified under $z\mapsto z^d$;
\item $\sum_{k=1}^s\big(\#(\ell_k\cap\partial\D)-1\big)=d-1$.
\end{enumerate}
A critical portrait $m$ is said to be a \emph{primitive major} if  point (1) in the definition above is strengthened to be that the elements of $m$ are pairwise disjoint.

\medskip
For $d = 2$, a primitive major is simply a diameter of the circle. Thus, each primitive major is parameterized by an angle $\theta \in \mathbb{R}/\mathbb{Z}$,
and the parameter space PM$(2)$ is homeomorphic to a circle.

\medskip
For $d = 3$, a cubic polynomial has either one critical point of multiplicity $2$, or two critical points of multiplicity $1$. Hence, there are two types of primitive majors: either an ideal triangle, or two leaves. For each pair $(a, b) \in S^1 \times S^1$ one can associate the pair of leaves $\{(a, a+1/3), (b, b+1/3)\}$, and
since the leaves cannot cross each other, one gets the strip
$$S := \left\{ (a, b) \in S^1 \times S^1 \ : \ a + \frac{1}{3} \leq b \leq a + \frac{2}{3} \right\}$$
which is the parameter space displayed in Figure \ref{F:global}. There are two particularly important slices:
\begin{enumerate}
\item The \emph{unicritical} slice $S_1$. It corresponds to the family $f(z) = z^3 + c$: each polynomial has only one critical point.
Combinatorially, it is represented by the slice $b = a + 1/3$ and the primitive majors are of the form $\{(a, a + 1/3, a + 2/3) \}$. See Figure \ref{F:uni}
for a graph of the restriction of the core entropy to this slice.

\item The \emph{symmetric} slice $S_2$. It corresponds to the family $f(z) = z^3 + c z$, and each polynomial is an odd function. Combinatorially,
it is represented by $b = a + 1/2$, so the associated primitive majors are $\{ (a, a + 1/3), (a + 1/2, a + 5/6) \}$. See Figure \ref{F:sym}.
\end{enumerate}

In order to account for the symmetries, note that $a$ and $b$ are interchangeable, so one can restrict to $a + \frac{1}{3} \leq b \leq a +\frac{1}{2}$ getting the
annulus
$$S':= \left\{ (a, b) \in S^1 \times S^1 \ : \ a + \frac{1}{3} \leq b \leq a + \frac{1}{2} \right\}.$$
The annulus $S'$ has two boundary components, corresponding to the two slices $S_1$ and $S_2$. One sees by the above discussion that the slice $S_1$
is periodic of period $3$, because the pairs $(a, a+1/3)$, $(a+1/3, a+2/3)$ and $(a+2/3, a)$ yield the same primitive major.
Moreover, the pairs $(a, a+1/2)$ and $(a+1/2, a)$ also yield the same primitive major, hence the slice $S_2$ is periodic of period $1/2$.

Thus, the parameter space PM$(3)$ is homeomorphic to the quotient of the annulus where one of the boundary components wraps around $2$ times, and the other boundary circle wraps around $3$ times. In formulas:
$${\rm PM}(3) = \left\{ (a, b) \in S^1 \times S^1 \ : \ a + \frac{1}{3} \leq b \leq a + \frac{1}{2} \right\}/\sim$$
where $(a, a+1/3) \sim (a+1/3, a+2/3) \sim (a+2/3, a)$ and $(a, a+1/2) \sim (a+1/2, a)$. The resulting space is not quite a manifold, as
a neighbourhood of the unicritical locus $S_1$ contains three ``sheets" which come together.

\subsection{Stratification}  \label{S:strata}

It is clear from the above discussion that ${ \rm PM}(d)$ has a natural stratification based on the size of the components of the primitive major.
Namely, let $s_1, \dots, s_r$ be integers with $\sum_{i = 1}^r (s_i - 1) = d- 1$. Then one can define the stratum $\Pi(s_1, \dots, s_r)$ as the set of primitive majors in ${ \rm PM}(d)$ which have leaves of size $s_1, \dots, s_r$.
In the above discussion of the cubic locus, the unicritical locus is the stratum $\Pi(3)$, while the generic stratum is $\Pi(2, 2)$.
A natural question then becomes:

\textbf{Question. }What is the maximum of core entropy on each stratum? How many (and which) polynomials achieve the maximum?

In general, the global maximum on ${\rm PM}(d)$ equals $\log d$, while as an example in the unicritical locus $\Pi(3)$ the maximum is $\log 2$.

As for the quadratic case, many other questions about the core entropy are completely open in higher degree, and it would be of great interest to pursue them. For instance, the local maxima of the core entropy are expected to occur at dyadic angles, while a conjecture analogous to the one in the quadratic case (formulated in \cite{Ti1} and proven in \cite{DS}) on the maximum of entropy on the wakes needs to be made precise. Moreover, we
believe the local H\"older exponent of the entropy function to be related to the value of the entropy (see \cite{Ti3} for real quadratic polynomials), and we expect self-similarity features in the graph of the entropy at preperiodic parameters (see \cite{Ju} for discussion of the quadratic case).

\subsection{Structure of the paper}

We start in Section 2 by reviewing the techniques in graph theory needed to define the core entropy through the \emph{spectral determinant}.
In Section 3, we define the combinatorial parameter space using primitive majors, and recall Thurston's entropy algorithm to compute the entropy.
Then (Section 4) we define the infinite graphs, called \emph{wedges}, which we use to encode the combinatorial dynamics in the space of postcritical arcs.
In order to study the limits of wedges as parameters vary, we define in Section 5 the concept of \emph{weakly periodic} labeled wedge.
In Section 7, we prove that any limit of sequence of wedges which correspond to a sequence of convergent parameters actually yields the same entropy.
Finally, in Section 8 we use this to establish the first main result, namely the continuity in the combinatorial parameter space (Theorem \ref{theorem:main}).

In the second part of the paper (Section 9) we transfer this combinatorial information to the analytic parameter space: there, we establish that as the coefficients of the polynomials converge,
then the critical markings also converge in a suitable way. This is achieved by showing continuity properties of the landing points of certain rays, and
introducing a generalization of the concept of critical marking (which we call \emph{weak critical marking}) which captures the marking of a limit
of postcritically finite polynomials. Using these tools we prove the second main result (Theorem \ref{theorem:main2}).

\subsection{Acknowledgements}

This paper is dedicated to the memory of Tan Lei (1963-2016).
We will always be grateful for her teachings, passion for the subject and encouragement.

Moreover, we wish to thank C. McMullen, J. Milnor, R. Perez, D. Thurston and M. Yampolsky for useful discussions.
G.T. is partially supported by NSERC and the Alfred P. Sloan Foundation.

\section{Growth rates of graphs of bounded cycles}

We start with some background material, following \cite[Sections 2,3]{Ti}.

\subsection{Graphs of bounded cycles}
In the following, by \emph{graph} we mean a directed graph, i.e., a set $V(\G)$ of
\emph{vertices} (which will be finite or countable) and a set $E(\G)$ of \emph{edges}, such that
each edge $e$ has a well-defined source $s(e)\in V$ and a target $t(e)\in V$ (thus,
we allow edges between a vertex and itself, and multiple edges between two
vertices). Given a vertex $v$, the set $Out(v)$ of its outgoing edges is the set of
edges with source $v$. The outgoing degree of $v$ is the cardinality of $Out(v)$;
a graph has \emph{bounded outgoing degree} if there is a uniform upper bound $d \geq 1$ on
the outgoing degree of all its vertices.

A path in the graph based at a vertex $v$ is a sequence $(e_1,\ldots,e_n)$ of edges such that $s(e_1)=v$ and $t(e_i)=s(e_{i+1})$ for $1\leq i\leq n-1$. The length of the path is the number $n$ of edges, and the set of vertices $\{s(e_1),\ldots,s(e_{n})\}\cup\{t(e_n)\}$ visited by the path is called its \emph{support}. Similarly, a closed path based at $v$ is a path $(e_1,\ldots,e_n)$ such that $t(e_n)=s(e_1)$. Note that in this definition closed
paths with different starting vertices will be considered to be different.

A \emph{simple cycle} is a closed path which does not self intersect, modulo cyclical equivalence: that is, a simple cycle is a closed path $(e_1,\ldots,e_n)$ such that $s(e_i)\not=s(e_j)$ for $i\not=j$, and two such paths are
considered the same simple cycle if the edges are cyclically permuted, i.e., $(e_1,\ldots,e_n)$ and $(e_{k+1},\ldots,e_n,e_1,\ldots,e_k)$. designate the same simple cycle.
Finally, a multi-cycle is the union of finitely many simple cycles with pairwise
disjoint (vertex-)supports. The length of a multi-cycle is the sum of the lengths of its components.

We say a graph has \emph{bounded cycles} if it has bounded outgoing degree and for each integer $n\geq1$ it has at most finitely
many simple cycles of length $n$.

Note that, if $\G$ has bounded cycles, then for each n it has also a finite
number of closed paths of length $n$. We shall denote as
\[C(\G,n)\]
the number of closed paths of length $n$, and define the \emph{growth rate} $r(\G)$ as the exponential growth rate of the number of its closed
paths: that is,
\[r(\G):=\limsup_{n\to\infty}\sqrt[n]{C(\G,n)}.\]

\subsection{The spectral determinant} Let $\G$ be a graph with bounded cycles. Let $S(\G,n)$ denote the number of simple multi-cycles of length $n$ in $\G$, and
let us define
\[\sigma:=\limsup_{n\to\infty}\sqrt[n]{S(\G,n)}\]
its growth rate. Tiozzo \cite[Section 2.1]{Ti} defined  a formal power series, called the \emph{spectral determinant}, as
\begin{equation}\label{eq:spetral-determinant}
P(t):=\sum\limits_{\g\text{ multi-cycle}}(-1)^{C(\g)}t^{\ell(\g)},
\end{equation}
where $\ell(\g)$ denotes the length of the multi-cycle, while $C(\g)$ is the number
of connected components of $\G$, and proved that the inverse of the growth rate $r(\G)$ is the minimal zero of $P(z)$:

\begin{lemma}[\cite{Ti}, Theorem 2.3]\label{lem:spetral-determinant}
Suppose we have $\sigma\leq1$; then the formula \eqref{eq:spetral-determinant} defines a holomorphic function $P(z)$ in the unit disk $|z|<1$, and moreover the function $P(z)$ is non-zero in the disk $|z|<r(\G)^{-1}$; if $r(\G)>1$, we also have $P(r(\G)^{-1})=0$.
\end{lemma}

For a finite graph $\G$ with vertex set $V=\{v_1,\ldots,v_q\}$, its \emph{adjacency matrix} $A=(a_{ij})_{q\times q}$ is defined as
\[a_{ij}:=\#(v_i\to v_j),\]
the number of edges from $v_i$ to $v_j$. In this case, the spectral determinant $P(t)$ equals ${\rm det}(I-tA)$.

\begin{lemma}[\cite{Ti}, Lemma 2.5]\label{lem:adjacent-matrix}
If $\G$ is a finite graph, then its growth rate equals the largest real eigenvalue of its adjacency matrix.
\end{lemma}

\subsection{Weak cover of graphs}

Let $\G_1,\G_2$ be two graphs with bounded cycles. A \emph{graph map} from $\G_1$ to $\G_2$
is a map $\pi:V(\G_1)\to V(\G_2)$ on the vertex sets and a map on edges $\pi:E(\G_1)\to E(\G_2)$ which is compatible, in the sense that if the edge $e$ connects $v$ to $w$ in $\G_1$, then the edge $\pi(e)$ connects $\pi(v)$ to $\pi(w)$ in $\G_2$.
We shall usually denote such a map as $\pi:\G_1\to\G_2$.

A \emph{weak cover} of graphs is a graph map $\pi:\G_1\to\G_2$
such that:
\begin{itemize}
\item the map $\pi:V(\G_1)\to V(\G_2)$ is surjective;
\item the induced map $\pi:Out(v)\to Out(\pi(v))$ between outgoing edges is a bijection for each $v\in V(\G_1)$.
\end{itemize}
As a consequence of the definition of weak cover, you have the following facts:
\begin{lemma}[\cite{Ti}, Lemma 3.1 and 3.3]\label{lem:fact-of-weak-cover}
Let $\pi:\G_1\to \G_2$ be a weak cover of graphs with bounded cycles. Then we have the following:
\begin{enumerate}
\item The unique path lifting property:  given $v\in\G_1$
and $w=\pi(v)\in\G_2$, for every path $\g$ in $\G_2$
based at $w$ there is a unique path $e$ in  $\G_1$ based at $v$ such that $\pi(e)=\g$;
\item  Let $S$ be a finite set of vertices of $\G_1$, and suppose that every closed path in $\G_1$ of length $n$ passes through $S$. Then we have the estimate
\[C(\G_1,n)\leq n\cdot\# S\cdot C(\G_2,n).\]
\end{enumerate}
\end{lemma}

\subsubsection{Quotient graphs}
A general way to construct weak covers of graphs
is the following. Suppose we have an equivalence relation $\sim$ on the vertex
set $V$ of a graph with bounded cycles, and denote $V^Q$ the set of equivalence classes of
vertices. Such an equivalence relation is called \emph{edge-compatible} if whenever $v_1\sim v_2$, for any vertex $w$ the total number of edges from $v_1$
to the members
of the equivalence class of $w$ equals that from $v_2$
to the members of the equivalence class of $w$. When we have such an equivalence
relation, we can define a quotient graph $\G^Q$  with vertex set $V^Q$ . Namely, we
denote for each $v,w\in V$ the respective equivalence classes as $[v]$ and $[w]$,
and define the number of edges from $[v]$ to $[w]$ in the quotient graph to be
\[\#([v]\to[w])=\sum\limits_{u\in[w]}\#(v\to u).\]
By definition of edge-compatibility, the above sum does not depend on the
representative $v$ chosen inside the class $[v]$. Moreover, it is easy to see that
the quotient map
\[\pi:\G\to \G^Q\]
is a weak cover of graphs.

\section{Core entropy for primitive majors}

\subsection{Primitive majors and critical portraits}\label{sec:major}

A \emph{critical portrait} of degree $d$ is a collection
\[m=\{\ell_1,\ldots,\ell_s\}\]
of leaves and ideal polygons in $\ov{\D}$ fullfilling the following conditions:
\begin{enumerate}
\item any two distinct elements $\ell_k$ and $\ell_l$ either are disjoint or intersect at one point on $\partial\D$;
\item the vertices of each $\ell_k$ are identified under $z\mapsto z^d$;
\item $\sum_{k=1}^s\big(\#(\ell_k\cap\partial\D)-1\big)=d-1$.
\end{enumerate}

The elements $(\ell_i)$ of a critical portrait will be called (portrait) \emph{leaves} (even when their cardinality is $> 2$, hence they correspond to polygons).
The number $s$ is called the \emph{size} of the critical portrait.

A critical portrait $m$ is said to be a \emph{primitive major} if  point (1) in the definition above is strengthened to be that the elements of $m$ are pairwise disjoint. We remark that for a critical portrait $m$ of degree $d$, the set $\D\setminus(\bigcup_{k=1}^s\ell_k)$ has $d$ connected components, and each one takes a total arc length $1/d$ on the unit circle (see \cite[Lemma 4.2]{G}).
 \begin{figure}[htpb]
\centering
\includegraphics[width=10cm]{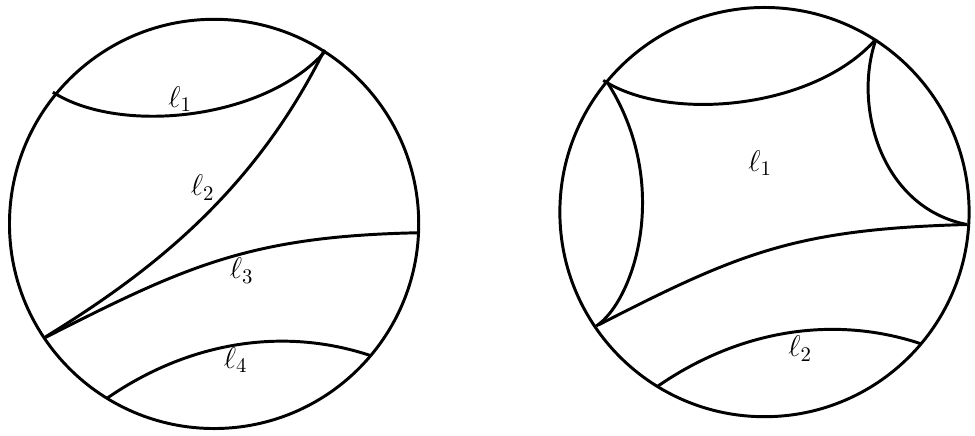}
\caption{Critical portraits of degree 5. The one on the right is a primitive major, while the one on the left is not a primitive major,
but it \emph{induces} the primitive major on the right. } \label{fig:major-portrait}
\end{figure}

In fact, each critical portrait induces a unique primitive major of the same degree. To see this, let $\xi=\{\ell_1,\ldots,\ell_s\}$ be a critical portrait of degree $d$. We define an equivalence relation on $\xi$ as the smallest equivalence relation such that if $\ell_i\cap \ell_j\not=\emptyset$, then $\ell_i$ and $\ell_j$ are equivalent. The portrait $\xi$ is therefore divided into the equivalence classes
$\QQQ_1,\ldots,\QQQ_t$. For each $i\in\{1,\ldots,t\}$, set
$$\Theta_i:=\text{the convex hull in $\ov{\D}$ of }\bigcup_{\ell_{j}\in\QQQ_i}(\ell_j\cap\partial\D).$$
The collection of sets
$\{\Theta_1,\ldots,\Theta_t\}$ is easily checked to be
 a degree $d$ primitive major, called the \emph{primitive major induced by $\xi$}.
For example, the primitive major in Figure \ref{fig:major-portrait} (right) is induced by the critical portrait on its left.

A critical portrait $\xi$ induces an equivalence relation on $\partial \mathbb{D}$, namely the smallest equivalence relation $\sim$
such that $x \sim y$ whenever $x$ and $y$ belong to the same leaf of $\xi$.
Two critical portraits are said to be \emph{equivalent} if they induce the same equivalence relation on $\partial \mathbb{D}$.
This is the same as saying that they induce the same primitive major.

\subsection{The topology in the space of primitive majors}

For $d\geq2$, we denote by ${\rm PM}(d)$ the space of all primitive majors of degree $d$. This space has a canonical metric md given by Thurston (see \cite[Part I, Section 3]{TG}) as follows.

\begin{figure}
 \centering
 \begin{tikzpicture}
 \node at (-4,1.5) {\includegraphics[width=4cm]{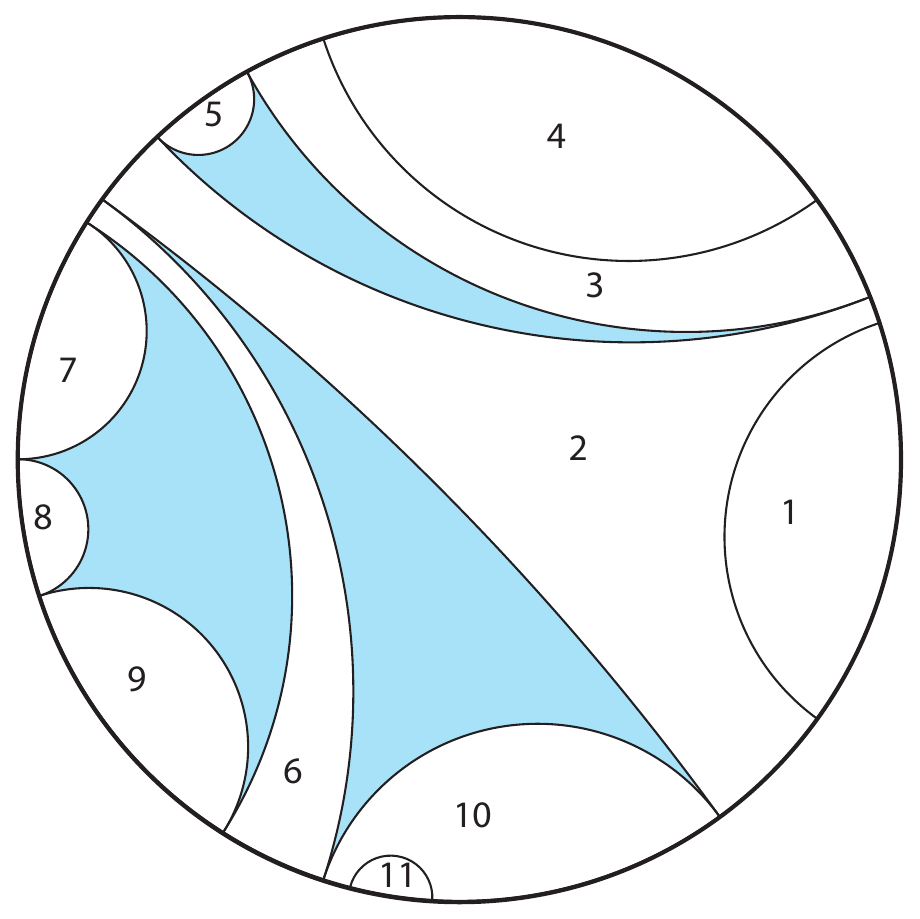}};
 \node at (3,1){\includegraphics[width=5.5cm]{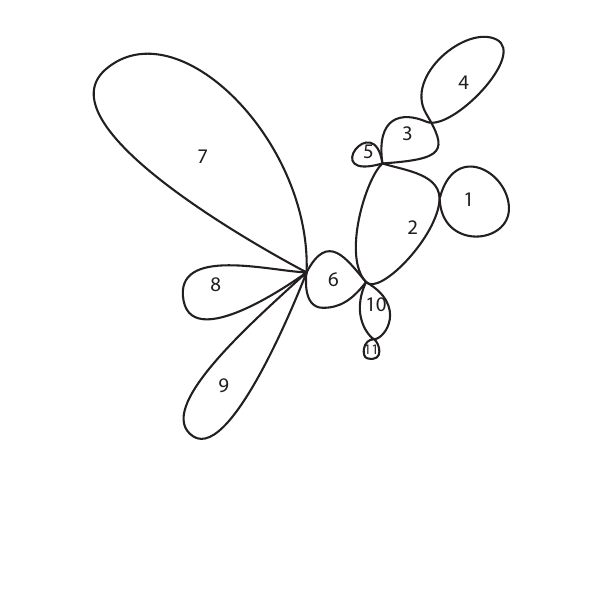}};
  \end{tikzpicture}
  \vspace{-1cm}
\caption{The deformation of the unit circle by shrinking a degree $7$ primitive major.}\label{shrink}
\end{figure}

A primitive major $m$ determines a quotient graph $\g(m)$ obtained from $\partial\D$ by identifying each element of $m$ to a point (see Figure \ref{shrink}). The path metric on $\partial \D$ determines a path metric on $\g(m)$.
Let met($m$) be the pseudometric on $\partial \D$ obtained as pullback of the path metric on $\gamma(m)$ under the projection $\partial \D \to \gamma(m)$;
then the
metric md on PM$(d)$ is defined as the sup difference of the (pseudo)metrics:
\[{\rm md}(m,m')=\sup_{x,y\in\partial\D}|{\rm met}(m)(x,y)-{\rm met}(m')(x,y)|.\]
We say that a sequence of majors $(m_n)_{n \geq 0}$ \emph{converges} to $m$ if the distance ${\rm md}(m_n, m)$ tends to zero.
\begin{example}\label{ex:convergence}
Let $d = 3$, and consider the following sequences of majors.
\begin{enumerate}
\item Set $m_n=\{\ell_1^n:=\{\frac{1}{n},\frac{1}{3}+\frac{1}{n}\},\ell_2^n:=\{-\frac{1}{n},\frac{2}{3}-\frac{1}{n}\}\}$. Then $(m_n)$
converges  to the primitive major $m=\{0,\frac{1}{3},\frac{2}{3}\}$.
\item Set $m_{2k}=\{\ell_1^{2k}:=\{\frac{1}{2k},\frac{1}{3}+\frac{1}{2k}\},\ell_2^{2k}:=\{-\frac{1}{2k},\frac{2}{3}-\frac{1}{2k}\}\}$ and $m_{2k-1}=\{\ell_1^{2k-1}:=\{\frac{1}{2k-1},\frac{1}{3}+\frac{1}{2k-1},\frac{2}{3}+\frac{1}{2k-1}\}\}$. Then $(m_n)$
also converges to the major $m$.
\end{enumerate}
\end{example}
From (2) of the example, we see that as  $(m_n)$ converges to $m$ in PM$(d)$, the size of $m_n$ may vary. Hence the sizes of their induced labeled wedges (see Section \ref{labeled-wedge}) may also vary. This will cause  difficulties in comparing the growth rates of the associated graphs.
To solve this problem, we partition the sequence $(m_n)$ into a finite number of convergent subsequences such that the majors in each subsequence have a common type.

For each critical portrait $\xi$, let $\xi^\cup$ denote the union of all portraits leaves of $\xi$, which is a compact subset of $\ov{\D}$.
A sequence of critical portraits $(\xi_n)$ is said to \emph{Hausdorff-converge} if the sequence of compact sets $(\xi^\cup_n)$ converges in the Hausdorff distance.

Note that, if the sequence of critical portraits $(\xi_n)$ Hausdorff-converges to a compact set $A$, then, for $n$ sufficiently large
we can label the elements of each $\xi_n$ by $\ell^{(n)}_1,\ldots,\ell^{(n)}_s$ such that for each $k\in\{1,\ldots,s\}$ the sequence $\ell_k^{(n)}$ of compact sets converges in the Hausdorff distance
to some compact set $\ell_k$, which is also either the closure of a leaf or of a polygon.
Moroever, one has $A=\bigcup_{k=1}^s\ell_k$. The set $\xi:=\{\ell_1,\ldots,\ell_s\}$ is called the \emph{limit} of $(\xi_n)$ in the sense of Hausdorff-convergence.

The following proposition shows the relation between convergence and Hausdorff-convergence for primitive majors.

\begin{proposition}\label{pro:limit-of-major}
Let $(\xi_n)$ be a sequence of critical portraits which Hausdorff-converge to $\xi$; we then have
\begin{enumerate}
\item the set $\xi$ is a critical portrait of degree $d$; and
\item if we let $m_n,n\geq1,$ be the primitive majors induced by $\xi_n$, and $m$  the primitive major induced by $\xi$, then the majors $(m_n)$ converge to $m$.
\end{enumerate}
Moreover, if the majors $(m_n)$ converge to $m$, the sequence $(m_n)$ can be partitioned into a finite number of Hausdorff convergent subsequences.
\end{proposition}

\begin{proof}
\begin{enumerate}
\item Since each $\ell_k$ is the Hausdorff limit of $\ell_k^{(n)}$, then we get that the vertices of $\ell_k^{(n)}$ converge to that of $\ell_k$, which implies that the vertices of $\ell_k$ are identified under $z\mapsto z^d$ and $\sum_{k=1}^s\big(\#(\ell_k\cap\partial\D)-1\big)=d-1$, and that each pair $\ell_k,\ell_l$ intersects at most on their boundary, which implies that the intersection is either one point on $\partial \D$ or a leaf. The latter case never happens because each component of $\D\setminus(\bigcup_{k=1}^s\ell^{(n)}_k)$ has total length $1/d$ on the unit circle. Then $\xi$ is a critical portrait of degree $d$.
\item Let $\zeta$ be any critical portrait of degree $d$, and $m_\zeta$ denote its induced primitive major. Note that the portrait $\zeta$ also induces a quotient graph $\g(\zeta)$ obtained from $\partial\D$ by identifying each element of $\zeta$ to a point, and it coincides with $\g(m_\zeta)$. Let met($\zeta$) be the pseudo-metric on the circle induced by the path-metric on $\g(\zeta)$. We then get ${\rm met}(\xi)={\rm met}(m_\zeta)$. For any $\epsilon>0$ and $x,y\in\partial\D$, as $(\xi_n)$ Hausdorff-converges to $\xi$, by the argument above, we have for $n$ large
    \[|{\rm met}(m_n)(x,y)-{\rm met}(m)(x,y)|=|{\rm met}(\xi_n)(x,y)-{\rm met}(\xi)(x,y)|<\epsilon.\]
    It follows that $(m_n)$ converges to $m$.
\item  We consider the accumulation points of the sequence $\{(m_n)^\cup\}$ in the Hausdorff topology.
Let $A$ be such an accumulation point and $(m_n')$ a subsequence with $(m'_n)^\cup\to A$ in the Hausdorff distance. It implies that the majors $m'_n$ Hausdorff-converge to $\xi$ with $\xi^\cup=A$. By assertions (1) and (2), such $\xi$ is a critical portrait which induces $m$. Note that the number of critical portraits that induce $m$ is finite, so the accumulation set of $\{(m_n)^\cup\}$ is finite.
Since the space of all compact subsets of a compact set is compact (in the Hausdorff topology), it follows that the sequence $(m_n)$ can be subdivided into a finite number of Hausdorff-convergent subsequences.
\end{enumerate}
\end{proof}

Note that the Hausdorff limit of primitive majors is not necessarily a primitive major: for example, the majors $(m_n)$ in Example \ref{ex:convergence} (1) Hausdorff-converge to the critical portrait $\xi=\{\{0,1/3\},\{0,2/3\}\}$. That is why we introduce the concept of critical portraits.

\subsection{Thurston's core entropy algorithm}\label{sec:algorithm}
We will describe here how Thurston's entropy algorithm works on rational primitive majors, see also \cite{G}.

By abuse of notation, we will identify a point of $\partial \D$ with its argument in $\T:=\R/\Z$. Then all angles in the circle are considered to be  mod 1, i.e. elements of $\T$.
The map $\tau:\T\to \T$  is defined by $\tau(\theta)=d\theta \mod \Z$.

Let $m=\{\ell_1,\ldots,\ell_s\}$ be a primitive major/critical portrait of degree $d$. Each $\ell_k$ is called a \emph{major/portrait} component. We set
\[x_k(i):= \tau^i(\ell_k\cap\T),\qquad i\geq1, k\in\{1,\ldots,s\}.\]
Note that the vertices of $\ell_k$ are identified by $\tau$, so all $x_k(i)$ are points in $\T$.
To describe Thurston's entropy algorithm, we consider the positions of $x_k(i)$ and $x_l(j)$ with respect to $m$.

 We say that the leaf $\ell$ \emph{separates}  two points $x_1$ and $x_2$ if $x_1$ and $x_2$ lie
in opposite connected components of $\T\setminus \ell$.

Given an ordered pair of points $(x_k(i), x_l(j))$ with $k,l\in\{1,\ldots,s\}$ and $i,j\geq1$, we say that their \emph{separation vector} (with respect to $m$) is
$(\alpha_1, \dots, \alpha_r)$ if the following are true:
\begin{enumerate}
\item each $\ell_{\alpha_i}$ belongs to $m$;
\item the leaf joining $x_k(i)$ and $x_l(j)$ successively crosses the leaves $\ell_{\alpha_1},\ldots,\ell_{\alpha_r}$ from $x_k(i)$ to $x_l(j)$ ;
\item no other element of $m$ separates $x_k(i)$ and $x_l(j)$.
\end{enumerate}
We say that $x_k(i)$ and $x_l(j)$ are \emph{not separated} if its separation vector is empty, and
they  are \emph{separated} otherwise.

The following fact will be used in the proof of Proposition \ref{pro:equi-portrait}.

\begin{lemma}\label{lem:separated}
 Let $\xi$ be a critical portrait, and $m$ the primitive major induced by $\xi$. Then two points of $\T$ are separated by an element $\Theta$ of $m$ if and only if they are separated by an element $\ell$ of $\xi$ which is contained in $\Theta$.
\end{lemma}

\begin{proof}
 The sufficiency is obvious. For the necessity,
  let $x,y\in \T$ be separated by an element $\Theta$ of $m$. If $\#\Theta=2$, then $\Theta$ is also an element of $\xi$, and the conclusion holds. Let $\#\Theta\geq3$. Then the leaf $\ov{xy}$ intersects two boundary leaves of $\Theta$, denoted by $\ov{\theta\eta}$ and $\ov{\theta'\eta'}$. Fixing $\theta$, there is one angle in $\{\theta',\eta'\}$, say $\theta'$, such that any connected subset of $\ov{\D}$ joining $\theta$ and $\theta'$ intersects $\ov{xy}$. Let $\ell$ and $\ell'$ be two elements of $\xi$ that contain $\theta$ and $\theta'$ respectively. By the construction of $m$ from $\xi$, there exist elements $\ell_{i_0}:=\ell,\ell_{i_1},\ldots,\ell_{i_t},\ell_{i_{t+1}}:=\ell'$ of $\xi$ contained in $\Theta$ such that $\ell_{i_k}\cap \ell_{i_{k+1}}\not=\emptyset$ for $k\in\{0,\ldots,t\}$. It follows that the connected set $\bigcup_{k=0}^{t+1}\ell_{i_k}$ joins $\theta$ and $\theta'$, and hence intersects $\ov{xy}$. Consequently, an element of $\xi$ among $\{\ell_{i_0},\ldots,\ell_{i_{t+1}}\}$ intersects $\ov{xy}$. Then the conclusion holds.
\end{proof}

A critical portrait $m=\{\ell_1,\ldots,\ell_s\}$ is said to be  \emph{rational} if all angles of
$\cup_{k=1}^s(\ell_k\cap\T)$ are rational numbers.

\textbf{The algorithm.} Let $\xi$ be a rational critical portrait. Then the set $\PPP(\xi):=\{x_k(i)\mid 1\leq k\leq s, i\geq1\}$ is finite. We define $O_\xi$ as the set of all unordered pairs
$\{x,y\}$ with $x\not=y\in \PPP(\xi)$ if $\#\PPP(\xi)\geq2$, and consisting of only $\{x,x\}$ if $\PPP(\xi)=\{x\}$. Then $O_\xi$ is finite but not empty.
The following is the procedure of Thurston's entropy algorithm acting on  $\xi$.
\begin{enumerate}
\item Let $\Sigma_{\xi}$ be the abstract linear space over $\R$ generated by the elements of $O_\xi$.

\item  Define a linear map $\mathcal{A}_{\xi}: \Sigma_{\xi}\longrightarrow \Sigma_{\xi}$ such that for any basis vector $\{x,y\}\in O_{\xi}$,
\begin{enumerate}
\item $\mathcal{A}_\xi(\{x,y\} )=0$ if $x,y$ belong to a common element $\ell_k$ of $\xi$;
\item $\mathcal{A}_\xi(\{x,y\} )=\{\tau(x),\tau(y)\}$ if $x,y$ are not separated by $\xi$ and do not belong to a common element of $\xi$; and
\item $\mathcal{A}_\xi(\{x,y\} )=\{\tau(x),  \tau(\ell_{\alpha_1})\} + \{\tau(\ell_{\alpha_1}),\tau(\ell_{\alpha_2})\}+\cdots+\{\tau(\ell_{\alpha_{r-1}}),\tau(\ell_{\alpha_r})\}+\{\tau(\ell_{\alpha_r}),\tau(y)\}$
if $x,y$ has separation vector $(\alpha_1,\ldots,\alpha_r)\not=\emptyset$.
 \end{enumerate}

\item  Denote by $A_\xi$ the matrix of $\mathcal{A}_\xi$ in the basis  $O_\xi$. It is a non-negative matrix. Compute its leading non-negative eigenvalue $\rho({\xi})$ (such an eigenvalue exists by the Perron-Frobenius theorem).  It is easy to see that $A_\xi$ is not nilpotent, therefore $\rho(\xi)\ge 1$.
\end{enumerate}

The output of Thurston's entropy algorithm is then
$$h(\xi) := \log \rho(\xi),$$
which we define as the \emph{core entropy} of the critical portrait $\xi$.

As proven by Gao Yan \cite{G}, the algorithm gives the correct value of the core entropy for \pf polynomials.
For a definition of weak critical marking, see Section \ref{section:critical-marking}.

\begin{theorem}[\cite{G}, Theorem 1.2]\label{entropy-algorithm-1}
Let $f$ be a \pf polynomial with  weak critical marking $\xi$. Then the core entropy $h(f)$ of $f$ is given by Thurston's algorithm, namely
$$h(f) = h(\xi).$$
\end{theorem}

\begin{figure}

\includegraphics[width = 0.5 \textwidth]{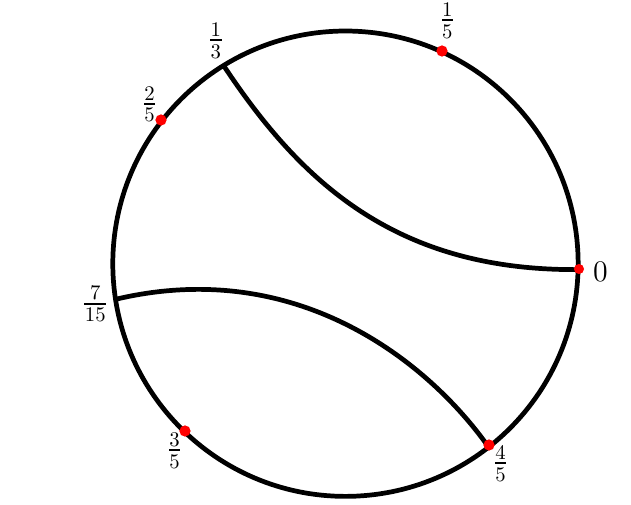}
\caption{An example of Thurston's entropy algorithm for a cubic polynomial.}
\label{F:algo}
\end{figure}

Let us conclude the section with an example of the algorithm, see also Figure \ref{F:algo}.

\begin{example}
Let $m=\{\ \{0,1/3\},\{7/15,4/5\}\ \}$. Then the set $\PPP(m)=\{0,1/5,2/5,3/5,4/5\}$
gives rise to an abstract linear space $\Sigma_m$ with  basis:
$$O_{m}= \Bigl\{\Bigl\{0,\frac{1}{5}\Bigr\},\Bigl\{0,\frac{2}{5}\Bigr\},\Bigl\{0,\frac{3}{5}\Bigr\},\Bigl\{0,\frac{4}{5}\Bigr\},\Bigl\{\frac{1}{5},\frac{2}{5}\Bigr\},
\Bigl\{\frac{1}{5},\frac{3}{5}\Bigr\},\Bigl\{\frac{1}{5},\frac{4}{5}\Bigr\},\Bigl\{\frac{2}{5},\frac{3}{5}\Bigr\},\Bigl\{\frac{2}{5},\frac{4}{5}\Bigr\},\Bigl\{\frac{3}{5},\frac{4}{5}\Bigr\}\Bigr\}$$
 The linear map $\AAA_{m}$ acts on the basis vectors as follows:
$$ \small{\begin{array}{l}
\Big\{0,\dfrac{1}{5}\Big\}\rightarrow\Big\{0,\dfrac{3}{5}\Big\},\quad \Big\{0,\dfrac{2}{5}\Big\}\rightarrow\Big\{0,\dfrac{1}{5}\Big\}, \quad \Big\{0,\dfrac{3}{5}\Big\}\rightarrow\Big\{0,\dfrac{2}{5}\Big\}+\Big\{\dfrac{2}{5},\dfrac{4}{5}\Big\},
\quad \Big\{0,\dfrac{4}{5}\Big\}\rightarrow\Big\{0,\dfrac{2}{5}\Big\},\\[15pt]
 \Big\{\dfrac{1}{5},\dfrac{2}{5}\Big\}\rightarrow\Big\{0,\dfrac{3}{5}\Big\}+\Big\{0,\dfrac{1}{5}\Big\}, \quad  \Big\{\dfrac{1}{5},\dfrac{3}{5}\Big\}\rightarrow\Big\{0,\dfrac{3}{5}\Big\}+\Big\{0,\dfrac{2}{5}\Big\}+\Big\{\dfrac{2}{5},\dfrac{4}{5}\Big\}, \quad \Big\{\dfrac{1}{5},\dfrac{4}{5}\Big\}\rightarrow\Big\{0,\dfrac{3}{5}\Big\}+\Big\{0,\dfrac{2}{5}\Big\},\\[15pt]
\Big\{\dfrac{2}{5},\dfrac{3}{5}\Big\}\rightarrow\Big\{\dfrac{1}{5},\dfrac{2}{5}\Big\}+\Big\{\dfrac{2}{5},\dfrac{4}{5}\Big\},\quad \Big\{\dfrac{2}{5},\dfrac{4}{5}\Big\}\rightarrow\Big\{\dfrac{1}{5},\dfrac{2}{5}\Big\}\quad  \Big\{\dfrac{3}{5},\dfrac{4}{5}\Big\}\rightarrow\Big\{\dfrac{4}{5},\dfrac{2}{5}\Big\}. \end{array}}$$
We compute  $h(m):=\log\rho(m)=1.395$.
\end{example}

\section{Labeled wedges and the associated graphs}\label{labeled-wedge}

We now turn to our general definition of core entropy for all critical portraits. In order to do so, we will generalize the transition matrix given by Thurston's algorithm to an infinite directed graph, which we call \emph{wedge}.

\subsection{Wedges} We will now introduce a combinatorial object, called a \emph{wedge}, to encode the dynamics on the set of postcritical arcs of a polynomial.
In fact, we will construct a graph whose vertices represent all possible arcs between the forward iterates of the critical points, and whose edges
represent the transitions between arcs as given by Thurston's algorithm.

Fix  an integer $d\geq2$. For each integer $1\leq s\leq d-1$, let us consider the set of pairs
$$\Sigma_s=\left\{ (y_k(i),y_l(j))  \ : \  k,l \in\{1,\ldots,s\}, 1 \leq i \leq j \right\},$$
which we call the {\it wedge} of size $s$. The pairs $(y_k(i),y_l(j))$ will be called \emph{vertices}, as they will become the vertices
of an infinite graph. In relation to the dynamics, the element $y_k(i)$ is meant to represent the $i^{th}$ iterate of the $k^{th}$ critical point,
while the vertex $(y_k(i), y_l(j))$ represents the arc between $y_k(i)$ and $y_l(j)$.

Moreover, given two elements $y_k(i)$ and $y_l(j)$ with $k, l \in \{1, \dots, s\}$ and $i,j \geq 1$, we denote as
$\{ y_k(i), y_l(j) \}$ the unique vertex in $\Sigma_s$ which represents an ordering of the pair consisting of $y_k(i)$ and $y_l(j)$ :
that is, $\{ y_k(i), y_l(j) \} = (y_k(i), y_l(j))$ if $i \leq j$, and $\{ y_k(i), y_l(j) \} = (y_l(j), y_k(i))$ otherwise.

\begin{figure}[h!]
$$\xymatrix{
& & & & \cdots \\
& & & (y_1(4),y_1(4)) & \cdots \\
& &  (y_1(3), y_1(3)) & (y_1(3),y_1(4))  &  \cdots \\
& (y_1(2),y_1(2)) & (y_1(2),y_1(3)) & (y_1(2),y_1(4)) & \cdots \\
(y_1(1),y_1(1)) & (y_1(1),y_1(2)) & (y_1(1), y_1(3)) & (y_1(1),y_1(4)) & \cdots \\
}$$
\caption{The wedge $\Sigma_1$: each vertex represent an arc between forward iterates of the critical point. For higher values of $s$, the wedge $\Sigma_s$ has several ``layers" given by the different values of $k, l$, which correspond to different critical points.}
\label{F:wedge}
\end{figure}

\subsection{Labeled wedges and graphs}
We then define a \emph{labeling} of the wedge of size $s$ as an assignment to each pair
$(y_k(i),y_l(j))$ of a label, which can be either $\emptyset$ or
$$(\alpha_1,\ldots,\alpha_r)$$
with $r\geq 1$ and $\alpha_1,\ldots,\alpha_r\in\{1,\ldots,s\}$ with the $\alpha_i$'s pairwise distinct.
We call the wedge of size $s$ with a labeling as above a {\it labeled wedge} of size $s$.
A vertex of $\Sigma_s$ is called \emph{non-separated} if its label is $\emptyset$, and \emph{separated} otherwise.

For each labeled wedge $\WWW$, we construct an \emph{associated graph} $\Gamma$ as follows. The vertices of $\Gamma$ are the elements of $\WWW$, and for each vertex $v= (y_k(i), y_l(j))$ of $\Sigma_s$, we determine the set of edges with source $v$ in the following way:
\begin{enumerate}
\item
if  $v$ is labeled $\emptyset$, then there is only one outgoing edge, namely
$$\{ y_k(i), y_l(j) \} \to \{ y_k(i+1), y_l(j+1) \};$$
such an edge will be called of \emph{upward} type.
\item
if the ordered pair $(y_k(i),y_l(j))$ is labeled  $(\alpha_1, \dots, \alpha_r)\not=\emptyset$, then
there are exactly $r+1$ edges going out of $v$, and precisely the following:
\begin{equation}\label{edge1}
\begin{array}{lll}
\{ y_k(i), y_l(j) \} & \to &
\{y_k(i+1), y_{\alpha_1}(1)\}\\
\{y_k(i), y_l(j)\} & \to &  \{y_{\alpha_1}(1), y_{\alpha_2}(1) \} \\
\{ y_k(i), y_l(j)\} & \to &  \{y_{\alpha_2}(1), y_{\alpha_3}(1) \} \\
& \dots &  \\
\{y_k(i), y_l(j) \} &  \to  & \{y_{\alpha_r}(1), y_l(j+1) \}.
\end{array}
\end{equation}
\end{enumerate}
By the definition of labeled wedge, the set of edges going out of $v$ is independent of the choice of the order of the elements of $v$.

\begin{figure}[h!]
$$\xymatrix{
& & & & \cdots \\
& & & (y_1(4),y_1(4)) & \cdots \\
& &  (y_1(3), y_1(3)) \ar[ur] & (y_1(3),y_1(4))  &  \cdots \\
& (y_1(2),y_1(2)) \ar[ur] & *+[F]{(y_1(2),y_1(3))}  \ar[dr] \ar[d] & (y_1(2),y_1(4)) & \cdots \\
(y_1(1),y_1(1)) \ar[ur] & *+[F]{(y_1(1),y_1(2))} \ar@(ul, ur) \ar[ur] & (y_1(1), y_1(3)) \ar[ur] & (y_1(1),y_1(4)) & \cdots \\
}$$
\caption{A labeling of the wedge $\Sigma_1$, and its associated graph. The boxed vertices have separation vector $\alpha_1$, representing the fact that
the (first) critical point lies on the associated arc, while the other vertices have empty separation vector.}
\end{figure}

For a vertex $v=\{y_k(i),y_l(j)\}\in\Sigma_s$, we call $\min\{i,j\}$ the \emph{height} of $v$, and $\max\{i,j\}$ the \emph{width} of $v$.
Among the first and last edges in the list \eqref{edge1}, the one whose target has smaller width will be called a \emph{backward} edge, and the other one is called a \emph{forward} edge.
All the other edges in \eqref{edge1} (the ones with targets of type $\{y_{\alpha}(1), y_{\beta}(1)\}$) will be called \emph{central} edges.

\subsection{The growth of labeled wedges}

\begin{proposition}\label{pro:bounded-cycle}
Let $\Gamma$ be the graph associated to a labeled wedge of size $s$. Then the following hold: 
\begin{enumerate}
\item each vertex along any closed path of length $n$ has height at most $n$;
\item each vertex along any closed path of length $n$ has width at most $2n$;
\item the number of simple multicycles of length $n$ is at most
$$(2nk)^{k + \sqrt{2kn}},$$
where $k := s^2$.
\end{enumerate}
\end{proposition}

\begin{proof}
\begin{enumerate}
\item Since the upward edges always increase the height of a vertex,
along each closed path there must be at least a backward, forward or a central edge.
Hence, since the target of a backward, forward
or central edge has height $1$, there must be at least one vertex of height $1$ along the closed path. Since
every edge increases the height of at most $1$, the claim follows.

\item Since the forward and upward edges always increase the width of a vertex,
along each closed path there must be at least a backward or central edge.
By the previous point, the source of such edge has height $\leq n$, hence its target has width $\leq n+1$.
The claim follows by the fact that each edge increases the width by at most $1$.
\item
Let $\gamma$ a simple multicycle of length $n$.
A vertex along the multicycle is called \emph{central} of type $(\alpha, \beta )$ if the edge of the multicycle
originating from it ends in the vertex $\{y_\alpha(1), y_\beta(1)\}$.
Moreover, a vertex is called \emph{backward} if it is separated and the backward edge originating from it belongs to the multicycle.

We first note that $\gamma$ is uniquely determined by the set of backward vertices, together with
the set of central vertices and their type.

Since the multicycle is simple, first note that along the multicycle there are at most $\frac{s(s-1)}{2} \leq k$ central vertices (at most one for each type).

We claim moreover that the number of backward vertices is at most $\sqrt{2 k n}$. In fact,
for each $h$ the number of backward vertices along the multicycle of height $h$ is at most $k := s^2$,
since the target of a backward edge whose source has height $h$ is of type $(y_\alpha(1), y_\beta(h+1))$,
and there are at most $s$ choices for $\alpha$ and $s$ choices for $\beta$.
Suppose now that the heights of the backward vertices along $\gamma$ are $h_1$, \dots, $h_r$.
Let us note that to each backward vertex of height $h_i$ there corresponds a segment of $\gamma$ of
length $h_i$, and all such segments are disjoint, so the total sum is $h_1 + \dots + h_r \leq n$.
Moreover, since for each value of $h$ there are at most $k$ values of $i$ such that $h_i = h$, we get the following estimate:
$$\frac{r^2}{2k} \leq \sum_{i =1}^r   \frac{i}{k} \leq \sum_{i = 1}^r h_i \leq n$$
which proves the upper bound on the number of backward vertices.

Finally, since there are at most $ 2 n k $ choices for each backward or central location
(at most one for each diagonal),
the total number of multicycles of length $n$ is at most $(2nk)^{k + \sqrt{2kn}}$, as required.
\end{enumerate}
\end{proof}

This proposition implies the growth rate of the number $S(\G,n)$ of simple multi-cycles is $\leq 1$. Then the following result follows directly from Lemma \ref{lem:spetral-determinant}.

\begin{corollary}
Let $\WWW$ be a labeled wedge, and let $\Gamma$ be its associated graph.
Then the graph $\Gamma$ has bounded cycles, and its spectral determinant $P(z)$ converges uniformly on compact subsets of the unit disk
$\mathbb{D} = \{ z \in \mathbb{C}, |z| < 1\}$, defining a holomorphic function $P: \mathbb{D} \to \mathbb{C}$. Moreover,
if the growth rate $r = r(\Gamma) > 1$, then the smallest real root of $P$ is $r^{-1}$. If $r  = 1$, then $P$ does not have any zeros in the unit disk.
\end{corollary}

We shall sometimes denote as $r(\WWW)$ the growth rate of the graph associated to the labeled wedge $\WWW$. We say that a sequence $(\WWW_n)_{n\geq1}$
of labeled
wedges of size $s$ \emph{converges} if for each finite set of vertices $S\subset \Sigma_s$ there exists $N$ such for each $n\geq N$ the labels of the elements of $S$ for $\WWW_n$ are the same.

\begin{lemma}[\cite{Ti}, Lemma~4.4]\label{lem:convergence-of-growth-rate}
If a sequence of labeled wedges $(\WWW_n)_{n\geq1}$ converges to $\WWW$, then the growth rate of $\WWW_n$ converges to that of $\WWW$.
\end{lemma}

\subsection{From critical portraits to labeled wedges}

We now see how to associate to each critical portrait a labeled wedge. Then, we will define the extension of the core entropy function as the growth
of the associated infinite graph.

Let $d \geq 2$, and let
$$\xi=\{ \ell_1, \dots, \ell_s \}$$
be a critical portrait of degree $d$. Recall that $x_k(i)=\tau^i(\ell_k)$ for each $k = 1, \dots, s$ and each $i \geq 1$.

The portrait $\xi$ induces a labeled wedge of size $s$ as follows:  for any vertex $\{y_k(i),y_l(j)\}$ of $\Sigma_s$, the ordered pair $(y_k(i),y_l(j))$ is labeled $(\alpha_1,\ldots,\alpha_r)$ ($r\geq0$) if  the ordered pair $(x_k(i),x_l(j))$  has the separation vector $(\alpha_1,\ldots,\alpha_r)$ with respect to $\xi$.

As an example, consider the critical portrait displayed in Figure \ref{F:algo}, with $d = 3$. We have
$$\xi = \left\{ \left(0, \frac{1}{3}\right), \left(\frac{7}{15}, \frac{4}{5}\right) \right\},$$
thus $\ell_1 = (0, 1/3)$ and $\ell_2 = (7/15, 4/5)$. As an example of labels, consider $x_2(2) = 3^2 \frac{4}{5} = \frac{1}{5} \mod 1$, and $x_2(3) = 3^3 \frac{4}{5} = \frac{3}{5} \mod 1$.
As one can see from the picture, the pair $(1/5, 3/5)$ is separated by the leaf $\ell_1 = (0, 1/3)$ and $\ell_2 = (7/15, 4/5)$, hence the
vertex $\{ y_2(2), y_2(3) \}$ has label $(\ell_1, \ell_2)$.
Then, the edges going out of this vertex are:
$$
\xymatrix{
&  \{y_2(3), y_1(1) \} \\
\{ y_2(2), y_2(3) \} \ar[r] \ar[ur] \ar[dr] &  \{y_1(1), y_2(1)\} \\
& \{y_2(1), y_2(4) \}
}
$$

We denote as $\WWW_\xi$ the labeled wedge induced by $\xi$, and as $\G_\xi$ its associated graph. The growth rate of $\G_\xi$ is simply denoted by $r(\xi)$.

\begin{definition}
Let $\xi$ be a critical portrait. Then the \emph{growth rate} $r(\xi)$ of $\xi$ is defined as the growth of the associated graph $\Gamma_\xi$.
\end{definition}

\subsection{Equivalence relation induced by a critical portrait}

So far we have constructed an infinite graph whose vertices represent all possible arcs joining forward iterates of the critical points.
However, iterates of postcritical angles may coincide. Thus, any critical portrait induces an equivalence relation on the circle, where
two pairs are defined to be equivalent if they represent the same pair of points on the circle. Let us see the details.

Let $\xi$ be a critical portrait.
We define an equivalence relation $\sim_\xi$ on the set $$\{y_k(i)\mid k\in\{1,\ldots,s\}, i\geq1\}$$  such that $y_k(i)\sim_\xi y_l(j)$  if $x_k(i)=x_l(j)$.
This means that the two forward iterates of the critical angles coincide.
This equivalence relation induces an equivalence relation, denoted by $\equiv_\xi$, on the vertices of the wedge $\Sigma_s$ such that 
$$\{y_{k_1}(i_1),y_{l_1}(j_1)\}\equiv_\xi\{y_{k_2}(i_2),y_{l_2}(j_2)\}$$
 if they are equivalent as a pair: that is, either
 $y_{k_1}(i_1)\sim_\xi y_{k_2}(i_2)\text{ and }y_{l_1}(j_1)\sim_\xi y_{l_2}(j_2)$, or
 $y_{k_1}(i_1)\sim_\xi y_{l_2}(j_2)\text{ and }y_{l_1}(j_1)\sim_\xi y_{k_2}(i_2).$

Finally, a vertex $v=\{y_k(i),y_l(j)\}\in \Sigma_s$ is called a \emph{diagonal vertex} with respect to $\xi$ if $y_k(i)\sim_\xi y_l(j)$: that is,
 the arc it represents is reduced to a single point.

\section{Weakly periodic labeled wedges}

 \begin{lemma}\label{pro:labeled-wedge-induced-by-portrait}
 Let $\xi$ be a critical portrait of degree $d\geq2$. Then the labeled wedge $\WWW_\xi$ satisfies
 \begin{enumerate}
 \item its diagonal vertices are all labeled $\emptyset$;
 \item if $v_1\equiv_\xi v_2\in\WWW_\xi$, then $v_1$ and $v_2$ have the same or the opposite label.
 \end{enumerate}
 \end{lemma}
\begin{proof}
It is easily checked by the definition of $\WWW_\xi$.
\end{proof}

As a generalization of the properties of $\WWW_\xi$ in Lemma \ref{pro:labeled-wedge-induced-by-portrait}, we get the concept of weakly periodic labeled wedge (of type $\xi$).

\begin{definition}
We call a labeled wedge $\WWW$  \emph{weakly periodic of type $\xi$} if the labels of its vertices
satisfy the following conditions.
\begin{enumerate}
\item

Suppose that the separation vector of the ordered pair $(x_k(i), x_l(j))$ is $(\alpha_1, \dots, \alpha_r)$.
Then the label in $\WWW$ of $(y_k(i), y_l(j))$ is of the form
\begin{equation}\label{eq:changed-label}
(\beta_1, \dots, \beta_t, \alpha_1, \dots, \alpha_r, \gamma_1, \dots, \gamma_u)
\end{equation}
where
$x_k(i) \in \beta_i$ for $i = 1, \dots, t$
and
$y_l(j) \in \gamma_i$ for $i = 1, \dots, u$.

Note that $t$ and $u$ may be zero, which shows that the condition is satisfied by the standard labeled wedge associated to $\xi$.

\item
\begin{enumerate}

\item
Moreover, if the label of the ordered pair $(y_k(i),y_{l}(j))$ is
$$(\beta_1, \dots, \beta_t, \alpha_1, \dots, \alpha_r, \gamma_1, \dots, \gamma_u)$$
and $x_l(j) = x_{l'}(j')$, then the label of the
ordered pair $(y_{k}(i),y_{l'}(j'))$  is
\begin{equation}\label{eq:change-3}
(\beta_1,\ldots,\beta_t ,\alpha_1,\ldots,\alpha_r, \gamma'_{1},\ldots,\gamma'_{u'})
\end{equation}
(i.e., the $\beta_i$ are the same).

\item
Similarly, if the label of the ordered pair $(y_k(i),y_{l}(j))$ is
$$(\beta_1, \dots, \beta_t, \alpha_1, \dots, \alpha_r, \gamma_1, \dots, \gamma_u)$$
and $x_k(i) = x_{k'}(i')$, then the label of the
ordered pair $(y_{k'}(i'),y_{l}(j))$  is
\begin{equation}\label{eq:change-3}
(\beta'_1,\ldots,\beta'_{t'} ,\alpha_1,\ldots,\alpha_r, \gamma_{1},\ldots,\gamma_{u}).
\end{equation}
(i.e., the $\g_i$ are the same).
\end{enumerate}
\end{enumerate}
\end{definition}

In the label \eqref{eq:changed-label}, we call the sub-vector
$(\beta_{1},\ldots,\beta_{t})$ the \emph{former-trivial labeled vector}, $(\gamma_{1},\ldots,\gamma_{u})$ the \emph{latter-trivial labeled vector}, and
$(\alpha_1,\ldots,\alpha_r)$ the \emph{essential labeled vector} of the ordered pair $(y_k(i),y_l(j))$.

Let $\G$ be the graph associated to a weakly periodic labeled wedge of type $\xi$. We denote $\G^{ND}$ the subgraph of $\G$ by taking as vertices all pairs which are non-diagonal, and as edges all the edges of $\G$ which do not have either as a source or target a diagonal pair.

\begin{lemma}\label{pro:compatible}
The equivalence relation $\equiv_\xi$ on $\G^{ND}$ is edge-compatible. Consequently, we get a quotient graph $\G^Q:=\G^{ND}/_{\equiv_\xi}$, and the quotient map
\[\pi:\G^{ND}\to \G^{Q}\] is a weak cover of graphs.
\end{lemma}
\begin{proof}
Let $v=\{y_k(i),y_l(j)\}$ and $v'=\{y_{k'}(i'),y_{l'}(j')\}$ be $\equiv_\xi$-equivalent. We assume that $y_k(i)\sim_\xi y_{k'}(i')$ and $y_{l}(j)\sim_\xi y_{l'}(j')$, and that the ordered pair $(y_k(i),y_l(j))$ has label as in eq. \eqref{eq:changed-label}.

 If the sub-vector $(\beta_{1},\ldots,\beta_{t})$ is not empty, we have
\[y_k(i+1)\sim_\xi y_{\beta_{1}}(1) \sim_\xi\cdots\sim_\xi y_{\beta_{t}}(1).\]
It follows that
\[\{y_k(i+1),y_{\beta_{1}}(1)\},\{y_{\beta_{1}}(1),y_{\beta_{2}}(1)\},\ldots\ldots,\{y_{\beta_{t-1}}(1),y_{\beta_{t}}(1)\}\]
are diagonal vertices, and the vertex $\{y_{\beta_{t}}(1),y_{\alpha_1}(1)\}$ is $\equiv_\xi$-equivalent to $\{y_k(i+1),y_{\alpha_1}(1)\}$. A similar argument holds for the sub-vector $(\g_{1},\ldots,\g_{u})$. Therefore, there are  at most $r+1$ edges in $\G^{ND}$ going out of $v$, and precisely all the ones
from following list which do not end in a diagonal vertex:
\begin{equation}\label{eq:all-edges}
\begin{array}{lll}
\{ y_k(i), y_l(j) \} & \overset{{e}_1}{\to} &
\left\{
  \begin{array}{ll}
    \{y_{\alpha_1}(1),y_k(i+1)\}, & \hbox{if $(\beta_{1},\ldots,\beta_{t})=\emptyset$;} \\
    \{y_{\alpha_1}(1),y_{\beta_{t}}(1)\}, & \hbox{otherwise.}
  \end{array}
\right.
\\
\{y_k(i), y_l(j)\} & \overset{{e}_2}{\to} &  \{y_{\alpha_1}(1), y_{\alpha_2}(1) \} \\
\{ y_k(i), y_l(j)\} & \overset{{e}_3}{\to} &  \{y_{\alpha_2}(1), y_{\alpha_3}(1) \} \\
& \dots &  \\
\{y_k(i), y_l(j) \} &  \overset{{e}_{r+1}}{\to}  &
\left\{
  \begin{array}{ll}
 \{y_{\alpha_r}(1), y_l(j+1) \}, & \hbox{if $(\g_{1},\ldots,\g_{u})=\emptyset$;} \\
    \{y_{\alpha_r}(1),y_{\g_{1}}(1)\}, & \hbox{otherwise.}
  \end{array}
\right.
\end{array}
\end{equation}

With a similar argument, we get that the edges in $\G^{ND}$ going out of $v'$ are precisely the non-diagonal ones among the following:
\begin{equation}
\begin{array}{lll}
\{ y_{k'}(i'), y_{l'}(j') \} & \overset{{e}_1'}{\to} &
\left\{
  \begin{array}{ll}
    \{y_{\alpha_1}(1),y_{k'}(i'+1)\}, & \hbox{if $(\beta_1',\ldots,\beta_{t'}')=\emptyset$;} \\
    \{y_{\alpha_1}(1),y_{\beta'_{t'}}(1)\}, & \hbox{otherwise.}
  \end{array}
\right.
\\
\{ y_{k'}(i'), y_{l'}(j') \} & \overset{{e}_2^\prime}{\to} &  \{y_{\alpha_1}(1), y_{\alpha_2}(1) \} \\
\{ y_{k'}(i'), y_{l'}(j') \} & \overset{{e}_3'}{\to} &  \{y_{\alpha_2}(1), y_{\alpha_3}(1) \} \\
& \dots &  \\
\{ y_{k'}(i'), y_{l'}(j') \} &  \overset{{e}_{r+1}'}{\to}  &
\left\{
  \begin{array}{ll}
 \{y_{\alpha_r}(1), y_{l'}(j'+1) \}, & \hbox{if $(\g'_1,\ldots,\g'_{u'})=\emptyset$;} \\
    \{y_{\alpha_r}(1),y_{\g'_1}(1)\}, & \hbox{otherwise.}
  \end{array}
\right.
\end{array}
\end{equation}

Note that, in any case,  the target of each $e_t$ ($1\leq t\leq r+1$) is $\equiv_\xi$-equivalent to the target of $e_t'$. It implies immediately  that the equivalence relation $\equiv_\xi$ is edge compatible.
\end{proof}

Note that by construction the quotient graph $\Gamma^Q$ can be also defined as follows.
Consider the postcritical set $\mathcal{P}$ on the circle:
$$\mathcal{P}=\PPP(\xi) := \{ \tau^i(\ell_k) \ : \ i \geq 1, 1 \leq k \leq s \}$$
Then take the set
$$\mathcal{A}=\AAA(\xi) := \{ \{x, y\} \in \mathcal{P} \times \mathcal{P} \ :  \ x \neq y \}$$
of non-degenerate pairs of postcritical points (the label $\mathcal{A}$ is because one thinks of it as the set of \emph{arcs} between postcritical points,
identifying an arc with its endpoints).
The set of vertices of $\Gamma^Q$ is precisely $\mathcal{P}$, while the set of edges is given by the dynamics.

\begin{proposition} \label{pro:equi-portrait}
Let $\xi_1$ and $\xi_2$ be two critical portraits, and $\Gamma_1$, $\Gamma_2$ be two weakly periodic labeled wedges of type, respectively, $\xi_1$ and $\xi_2$. If $\xi_1$ is
equivalent to $\xi_2$ (see Section \ref{sec:major}), then the quotient graphs $\Gamma_1^Q$ and $\Gamma_2^Q$ are isomorphic.
\end{proposition}

\begin{proof}
Since any critical portrait is equivalent to exactly one primitive major, there exists a
primitive major $m$ which is equivalent to both $\xi_1$ and $\xi_2$. Thus, it is enough to prove the statement when one of the two critical portraits, say $\xi_1$, is a primitive major.

If two leaves intersect on the boundary, then they have the same image under $\tau$. Hence the set of images
$\mathcal{P}  = \{ \tau^i(\ell_k) \ : \ i \geq 1, 1 \leq k \leq s \}$
is the same for $\xi_1$ and $\xi_2$. Thus, the two graphs $\Gamma_i^Q$ have the same vertex set.

In order to check the edges, let us now consider a vertex $v = (x_k(i), x_l(j))$ of $\Gamma_1^Q$ (and also of $\Gamma_2^Q$, as seen above),
 and let us suppose that the separation vector of the two points $x = x_k(i)$ and $y = x_l(j)$ on the circle equals $(\alpha_1, \dots, \alpha_r)$.
By definition of weakly periodic, the label of $v$ in $\Gamma^Q_2$ equals
$(\beta_1, \dots, \beta_t, \alpha_1, \dots, \alpha_r, \gamma_1, \dots, \gamma_u)$ for some choice of $\beta_i$ and $\gamma_i$.
Then note that, since $x$ belongs to all the $\beta_i$, the arcs $(\tau(\beta_i), \tau(\beta_{i+1}))$ for $i = 1, \dots, t-1$
are all degenerate. So are the arcs $(\tau(\gamma_i), \tau(\gamma_{i+1}))$ for $i = 1, \dots, u-1$ since $y$ belongs to all the $\gamma_i$.
Thus, the outgoing edges from $v$ are the non-degenerate arcs among
$(\tau(x), \tau(\alpha_1)), (\tau(\alpha_1), \tau(\alpha_2)), \dots, (\tau(\alpha_{r-1}), \tau(\alpha_r)), (\tau(\alpha_r), \tau(y))$.
Now, by definition of the equivalence relation there exists equivalence classes $\Theta_1, \dots, \Theta_w$
such that
$\alpha_1, \dots, \alpha_{i_1}$ belongs to $\Theta_1$,
$\alpha_{i_1+1}, \dots, \alpha_{i_2}$ belongs to $\Theta_2$, etc.
Then we note that the arcs $(\tau(\alpha_1), \tau(\alpha_2))$, up to $(\tau(\alpha_{i_1-1}), \tau(\alpha_{i_1}))$ are also degenerate, and so on,
hence the outgoing edges from $(x,y)$ are the non-degenerate arcs among
$$(\tau(x), \tau(\alpha_1)), (\tau(\alpha_1), \tau(\alpha_2)), \dots, (\tau(\alpha_{r-1}), \tau(\alpha_r)), (\tau(\alpha_r), \tau(y))$$
This is by definition the list of outgoing edges from $(x, y)$ in $\Gamma_1^Q$, proving the claim.
\end{proof}

\section{The comparison of growth rates of weakly periodic labeled wedges}

Throughout this section, we always assume that $$\xi=\{\ell_1,\ldots,\ell_s\}.$$

Let $\WWW$ be a weakly periodic labeled wedge of type $\xi$, and $\G$ its associated graph. The notations  $\G^Q$ and $\G^{ND}$ follow Proposition \ref{pro:compatible}.

\begin{proposition}\label{lem:common-growth-rate}
Let $\WWW$ be a weakly periodic labeled wedge of type $\xi$, and $\G$ its associated graph.
Then the growth rates of $\G^{ND}$ and $\G^Q$ are equal.
\end{proposition}

\begin{proof}
Let $S$ denote the set of vertices in $\Sigma_s$ which have widths and heights at most $2n$. By Proposition \ref{pro:bounded-cycle}, each closed path in $\G^{ND}$ of length $n$ passes through $S$. Applying (2) of Lemma \ref{lem:fact-of-weak-cover}, we get the estimate
\[r(\G^{ND})\leq r(\G^Q).\]
We then need to show $r(\G^Q)\leq r(\G^{ND})$.

Let $\g=(e_1,\ldots,e_n)$ be a closed path in $\G^Q$ with  ${v}_0:=s({e}_1)$ and ${v}_t:=t({e}_t)$ for each $t\in\{1,\ldots,n\}$.
Each $v_t$ represents an arc, hence it has two endpoints.
By induction on $t$, we will declare certain endpoints of $v_t$ as \emph{marked}, according to the following rule:

\begin{enumerate}
\item by definition, both endpoints of the arc $v_0$ are marked.

\item recursively, an endpoint of $v_{t+1}$ is marked if it is the image of a marked endpoint of $v_t$, in the following sense.
Let $v_t = (x_k(i), x_l(j))$, and suppose the separation vector is $(\alpha_1, \dots, \alpha_r)$.
Then, if the endpoint $x_k(i)$ is marked, we also mark the endpoint $x_k(i+1)$ of the arc  $(x_k(i+1), x_{\alpha_1}(1))$.
Similarly, if $x_l(j)$ is marked, then we mark the endpoint $x_l(j+1)$ of the arc $(x_{\alpha_r}(1), x_l(j+1))$.
All other endpoints of $v_{t+1}$ are not marked.
\end{enumerate}

We then say that a vertex $v_t$ is marked if at least one of its endpoints is marked.
A path $\gamma$ with vertices $v_0, \dots, v_n$ is \emph{peripheral} if all its vertices are marked, and \emph{non-peripheral} otherwise.

Note that by construction, out of all the edges going out of $v_t$, at most two are marked.
As a consequence, for any $v_0$ and any $n$, there are at most two peripheral paths of length $n$ which start at $v_0$.

\textbf{Claim.}
If $\g$ is peripheral, then there is a vertex $\wt{v}_0\in\pi^{-1}(v_0)$ which has width and height at most $C n$, where $C$ is a constant which depends
only on $\xi$.

Let us pick a vertex $w_0 = \{y_k(i),y_l(j)\}\in\pi^{-1}(v_0)$, and consider the lift $\wt{\g}=(\wt{e}_1, \ldots, \wt{e}_n)$ of $\g$ based at $w_0$.
Set $w_t:=t(\wt{e}_t)$ for each $t\in\{1,\ldots,n\}$. The claim will be checked by cases.
\begin{enumerate}
\item Each vertex $w_0,\ldots, w_{n-1}$ is labeled $\emptyset$. Then $w_n=\{y_k(i+n),y_l(j+n)\}$. Note that $w_0\equiv_\xi w_n$, so we get either that $x_k(i)=x_k(i+n)$ and $x_l(j)=x_l(j+n)$ or that $x_k(i)=x_l(j+n)$ and $x_l(j)=x_k(i+n)$. In both cases $\ell_k$ and $\ell_l$ are eventually periodic. It follows that there are a constant $C_1>0$ and $i_0,j_0<C_1$ such that $x_k(i_0)=x_k(i)$ and $x_l(j_0)=x_l(j)$. The point $\wt{v}_0:=\{y_k(i_0),y_l(j_0)\}$ thus satisfies the requirements.
\item There is a separated vertex among $w_0,\ldots,w_n$. To better show the argument, let us first assume that $w_0$ is separated. Then $w_1$ has height $1$.

If there is  a central or backward edge among $\wt{e}_2,\ldots, \wt{e}_n$, then the width and height of $w_n$ are both less than $n$, and $\wt{v}_0:=w_n$ satisfies the requirement.

Otherwise, we get that $w_1$ equals $\{y_\alpha(1),y_{k}(i+1)\}$ or $\{y_\beta(1),y_l(j+1)\}$, and the edges $\wt{e}_2,\ldots,\wt{e}_n$ are either forward or upward.
By symmetry, we can assume $w_1=\{y_\alpha(1),y_{k}(i+1)\}$. It follows that $w_n=\{y_{\alpha'}(p), y_k(i+n)\}$ with $p\leq n$ and $w_n\equiv_\xi w_0$. If $y_k(i)\sim_\xi y_{\alpha'}(p)$ and $y_l(j)\sim_\xi y_k(i+n)$,  we set
$$\wt{v}_0:=\{y_{\alpha'}(p),y_{\alpha'}(p+n)\},$$
 which is $\equiv_\xi$-equivalent to $w_0$, and has width and height at most $2n$. If $y_k(i)\sim_\xi y_k(i+n)$ and $y_l(j)\sim_\xi y_{\alpha'}(p)$, then $\ell_k$ is eventually periodic. There is hence an integer $i_0$ less than a constant $C_2$ such that $y_k(i_0)\sim_\xi y_k(i)$. The vertex
 $$\wt{v}_0=\{y_{\alpha'}(p),y_k(i_0)\}$$
 is what we want, with width and height at most $C_2n$.

In the general case, let $w_k$ be the first separated vertex among $w_0, \dots, w_n$. Then by the previous argument there exists a vertex $\wt{v}_k$ which
projects to $v_k$ and has height and width $\leq \max\{C_1, C_2\} n$. By lifting the path $(e_{k+1}, \dots, e_{n})$ starting from $w_k$ one gets a vertex $\wt{v}_0$ which projects to $v_0$ and with height and width bounded above by $\max\{C_1, C_2\} n + n$, as required.

\end{enumerate}

Now, let us note that the number of vertices of the wedge with both width and height bounded by $C n$ is at most $(s C n)^2$, hence
the number of projections to $\Gamma_Q$ of such vertices is also bounded above by $(s C n)^2$. Finally, as we previously observed the number of
peripheral paths of length $n$ starting at a given vertex $v_0$ is at most $2$, hence we get the estimate

\begin{equation} \label{E:non-per}
\Big\{
\begin{array}{ll}
\text{peripheral closed paths}\\
\text{ of length $n$}
\end{array}
\Big\}\leq 2 s^2 C^2 n^2.
\end{equation}

\textbf{Claim.} If $\g$ is a non-peripheral closed path in $\G^Q$, then there exists a closed path $\wt{\g}\subset \G^{ND}$ of length $n$
which projects in $\Gamma^Q$ to a cyclic permutation of $\gamma$. 

\begin{proof}[Proof of the Claim]
By cyclic permutation of $\gamma = (e_1, \dots, e_n)$, we mean a path of the form $(e_{k}, \dots, e_n, e_1, \dots, e_{k-1})$ for some $k$.
By definition of non-peripheral, the exists $n_1 \leq n$ the least index for which at least one of the endpoints of $v_{n_1}$ is not marked,
and $n_2 \in [n_1, n]$ the least index for which none of the endpoints of $v_{n_2}$ is marked.
Let us choose now a vertex $\widetilde{v}_0$ of $\Gamma^{ND}$ which projects to $v_{0}$, and let us lift $\gamma$ starting from there.
Thus we get a sequence $\widetilde{v}_0, \widetilde{v}_1, \dots, \widetilde{v}_n$ of vertices, and it is not necessarily true that $\widetilde{v}_0 = \widetilde{v}_n$. Let us now keep lifting $\gamma$ starting from $\widetilde{v}_n$, obtaining a sequence $\widetilde{v}_{n+1}, \dots, \widetilde{v}_{2n}$
 which also projects to $\gamma$.

Now, by definition of weakly periodic labeled wedge, the two vertices $\widetilde{v}_{n_1}$ and $\widetilde{v}_{n+n_1}$ have a common endpoint.
Then, by applying successively conditions (2) (a)-(b) of the definition, one gets that the same is true for all the pairs $\widetilde{v}_{t}$ and
$\widetilde{v}_{t+n}$ with $n_1 \leq t < n_2$. Finally, this implies that $\widetilde{v}_{n_2} = \widetilde{v}_{n+n_2}$.

Thus, by lifting the path $v_{n_2}, \dots, v_{n}, v_{1},\dots, v_{n_2-1}, v_{n_2}$ in $\Gamma^Q$ one gets a closed path in $\G^{ND}$, as required.
\end{proof}

By equation \eqref{E:non-per} and the previous claim, we have the estimates
\begin{eqnarray*}
C(\G^{Q},n)&=&\#\Big\{
\begin{array}{ll}
\text{peripheral closed}\\
\text{paths of length $n$}
\end{array}
\Big\}
+\#\Big\{
\begin{array}{ll}
\text{non-peripheral closed}\\
\text{paths of length $n$}
\end{array}
\Big\}\\
&\leq& 2 s^2 C^2 n^2+n C(\G^{ND},n),
\end{eqnarray*}
from which follows
\[r(\G^{Q})\leq r(\G^{ND})\]
as required.
\end{proof}

\begin{lemma} \label{L:same-entro}
Let $\xi$ be a critical portrait, and $m$ its induced primitive major. Then the equation
$$r(\xi)=r(m)$$
holds.
\end{lemma}

\begin{proof}
Note that the labels of every diagonal of the labeled wedges $\mathcal{W}_\xi, \mathcal{W}_m$ are empty, so we have $r(\xi):=r(\G_\xi)=r(\G_\xi^{ND})$ and $r(m):=r(\G_m)=r(\G_m^{ND})$. It then follows from Propositions \ref{lem:common-growth-rate} and \ref{pro:equi-portrait} that
 \[r(\xi)=r(\G_\xi^{ND})\overset{{\rm Pro}.\ref{lem:common-growth-rate}}{=}r(\G^Q_\xi)\overset{{\rm Pro.}\ref{pro:equi-portrait}}{=}r(\G^Q_m)\overset{{\rm Pro}.\ref{lem:common-growth-rate}}{=}r(\G_m^{ND})=r(m),\]
 which proves the claim.
\end{proof}

\begin{lemma} \label{L:wplimit}
Let $(m_N)$
be a sequence of primitive majors which Hausdorff-converge to a critical portrait $\xi$, and suppose that the associated
sequence of labeled wedges $(\mathcal{W}_N)$ converges to some labeled wedge $\mathcal{W}$. Then $\mathcal{W}$ is weakly periodic of type $\xi$.
\end{lemma}

\begin{proof}
In order to check (1) of the definition of weakly periodic labeled wedge, let the ordered pair $(x_k(i),x_l(j))$ have the separation vector $(\alpha_1,\ldots,\alpha_r)$ with respect to $\xi$.
Note that if
$x_k(i)\not\in\ell_\alpha$, then for $N$ large the point $x^{(N)}_k(i)$ has the same position with respect to $\ell_{\alpha}^{(N)}$ as that of $x_k(i)$ with respect to $\ell_\alpha$.  Thus, the leaves $\alpha_1, \dots, \alpha_r$ must be part of the separation vector of $(x_k^{(N)}(i),x_l^{(N)}(j))$, and on the other hand
the only other leaves which are part of this separation vector must contain either $x_k(i)$ or $x_l(j)$.  Since $\mathcal{W}_N \to \mathcal{W}$,
the label of $(y_k(i), y_l(j))$ in $\mathcal{W}$ equals the separation vector for  $(x_k^{(N)}(i), x_l^{(N)}(j))$  for $N$ large, so
this argument proves property (1) in the definition of weakly periodic labeled wedge.

Let us now prove (2)(a).
Let $v'=\{y_{k}(i),y_{l'}(j')\}\in \Sigma_s$ be a vertex $\equiv_\xi$-equivalent to $v$, i.e. so that $x_l(j) = x_{l'}(j')$, and let
the separation vector of $(x_{k}(i),x_{l}(j))$ be $(\alpha_1, \dots, \alpha_r)$.
Then the separation vector of $(x_k^{(N)}(i), x_l^{(N)}(j))$ is of type $(\beta_1, \dots, \beta_t, \alpha_1, \dots, \alpha_r, \gamma_1, \dots, \gamma_u)$,
where the $\ell^{(N)}_{\beta_i}$ are precisely the leaves which separate $x_k^{(N)}(i)$ and $\ell_{\alpha_1}^{(N)}$.
For the same reason, the separation vector of $(x_k^{(N)}(i), x_{l'}^{(N)}(j'))$ is of type $(\beta_1, \dots, \beta_t, \alpha_1, \dots, \alpha_r, \gamma'_1, \dots, \gamma'_{u'})$ (note the $\beta_i$ are the same). Since $\mathcal{W}_N$ converges to $\mathcal{W}$, these are also the labels of, respectively,
 $(y_k(i), y_l(j))$ and $(y_k(i), y_{l'}(j))$ in $\mathcal{W}$, proving the claim. (2) (b) follows analogously.
\end{proof}

\section{The convergence of labeled wedges induced by primitive majors}

To prove the continuity of the growth rate $r(m)$ (about primitive majors $m$ in the metric md), we expect to apply Lemma \ref{lem:convergence-of-growth-rate}. For this purpose, we need to know when the labeled wedges $\WWW_{m_N}$ converge as the majors $m_N$ converge. Note that even if the $m_N$ converge in the Hausdorff topology,
 the labeled wedges $\WWW_{m_N}$ may not converge. For example, in the quadratic case, if $\theta$ is periodic, then the labeled wedges $\WWW_{\theta'}$ do not converge as $\theta'\to\theta$. However, we will show that it is true for a subsequence.

\begin{lemma} \label{L:subseq}
Let $s \geq 1$. Then any sequence $(\mathcal{W}_N)$ of labeled wedges of size $s$ has a convergent subsequence.
\end{lemma}

\begin{proof}
It follows by our choice of (weak!) topology on the space of labeled wedges.
Since any vertex of $\Sigma_s$ has finitely many possible labels,
for each finite set $S$ of vertices of $\Sigma_s$ there exists a subsequence $(\mathcal{W}_{N_k})$ of labeled wedges
such that all vertices of $S$ have the same label. The claim follows by picking an
exhaustion of $\Sigma_s$ by finite sets $(S_n)$ and applying the usual diagonalization argument.
\end{proof}

\begin{lemma}\label{lem:limit-labeled wedge}
Let $(m_N)$ be a sequence of primitive majors which converges to a critical portrait $\xi$ in the Hausdorff topology, and so that
the associated labeled wedges $\mathcal{W}_{m_N}$ converge to a wedge $\mathcal{W}_\infty$, with associated infinite graph $\Gamma_\infty$. Then
$$r(\Gamma_\infty) = r(\Gamma_\infty^{ND}).$$
\end{lemma}

\begin{proof}

Let us denote as $\ell_1, \dots, \ell_s$ the leaves of $\xi$, and denote as $x_k(i)$ the point on the circle $x_k(i) = \tau^i(\ell_k)$. Recall that a vertex
$v = \{ y_k(i), y_l(j) \}$ of $\Gamma_\infty$ is diagonal with respect to $\xi$ if $x_k(i) = x_l(j)$.

For the graph $\G$ associated to a labeled wedge and any integer $n\geq1$, we denote $\G^n$ the finite subgraph of $\G$ such that $V(\G^n)$
is the set of vertices of $\G$ with width and height at most $n$, and $E(\G^n)$ is the set of all edges of $\G$ with both sources and targets in $V(\G^n)$.

For a primitive major $m$, we will denote as $\WWW_m$ the associated labeled wedge, and as $\Gamma_m$ its associated infinite graph.
Let us now fix $n \geq 1$.
Since $\WWW_{m_N} \to \WWW_\infty$, then we can choose $m = m_N$ sufficiently close to $\xi$ so that each vertex $v\in\Sigma_s$ with width and height at most $2n$ has a common label in $\WWW_m$ and $\WWW_\infty$. It follows that all graphs $\G_m^{2n}$ coincide with $\G_\infty^{2n}$.  Moreover, note that by Proposition \ref{pro:bounded-cycle} every closed path of length $n$ in $\Gamma_\infty$ actually lives in $\Gamma_\infty^{2n}$, which is also equal to $\Gamma_m^{2n}$.

To prove $r(\G_\infty)=r(\G^{ND}_\infty)$, we shall check the estimate
\[\#\Big\{
\begin{array}{ll}
\text{closed paths of length $n$ in $\G^{2n}_\infty$,}\\
\text{containing $\xi$-diagonal vertices}
\end{array}
\Big\}\leq 4s^2n^2\]
which by the above discussion is equivalent to the estimate
\[\#\Big\{
\begin{array}{ll}
\text{closed paths of length $n$ in $\G^{2n}_m$,}\\
\text{containing $\xi$-diagonal vertices}
\end{array}
\Big\}\leq 4s^2n^2.\]
This result will follow from the following fact:

$(\star)$ For each $n$ and each diagonal vertex $v_0$ of $\Gamma_m$, there exists at most one closed path of length $n$ based at $v_0$ in $\Gamma_m$.

First note that the fact implies the claim,
as the height and width of $v_0$ are bounded above by $2n$, yielding the estimate
\[\#\Big\{
\begin{array}{ll}
\text{closed paths of length $n$ in $\G^{2n}_m$,}\\
\text{containing diagonal vertices}
\end{array}
\Big\}\leq (2n\cdot s)^2=4s^2n^2\]
as required.

Let us now prove $(\star)$. In order to do so, for each major $m$ approximating $\xi$ let us denote as $x^m_k(i) := \tau^i(\ell^m_k)$ the iterate of the approximating leaf $\ell^m_k$.

Suppose that $v_0  = \{ y_{k_0}(i_0), y_{l_0}(j_0) \}$ is a $\xi-$diagonal vertex, and let
 $\theta_0 = x_{k_0}(i_0) = x_{l_0}(j_0)$. Let us choose an interval $I$ in the circle which contains $\theta_0$ in its interior, and such that the map $\tau^n : I \to \tau^n(I)$ is a homeomorphism.

Let us now choose a primitive major $m = m_N$ close enough to the limit $\xi$ so that  $x^m_{k_0}(i_0)$ and $x^m_{l_0}(j_0)$ belong to $I$.
Note that if $x^m_{k_0}(i_0)$ and $x^m_{l_0}(j_0)$ coincide, then the vertex $v_0$ is not separated in $\Gamma_m$, and so are all its descendants in the graph $\Gamma_m$:
thus, $v_0$ does not lie on any closed path. Thus, we can assume that the interval $[x^m_{k_0}(i_0), x^m_{l_0}(j_0)]$ is not a point.
	
For each vertex $v = \{ y_k(i), y_l(j) \}$ and each approximating major $m$, let us denote $J_v^m := [ x^m_k(i), x^m_l(j) ]$
the corresponding arc on the circle connecting the two iterates of the approximating major.

Suppose now that there is a path $v_0 \to v_1 \to \dots \to v_{n-1} \to v_n$ in $\Gamma_m$.
Note that by construction each interval $J^m_{v_{t+1}}$ is a subinterval of $\tau(J^m_{v_t})$, thus
$J^m_{v_n}$ is a subinterval of $L = \tau^n(J^m_{v_0})$.  Moreover, distinct paths yield disjoint subintervals.

If the path is closed ($v_0 = v_n)$, then the intervals $J^m_{v_0}$ and $J^m_{v_n}$ must coincide. However, as all intervals $J^m_{v_n}$
for different choices of paths are disjoint, there is at most one path for which   $J^m_{v_n}$ coincides with $J^m_{v_0}$.
Thus, there exists at most one
closed path of length $n$ based at $v_0$, proving $(\star)$.
\end{proof}

\section{The continuity of growth rate and core entropy}

In this part, we will show the continuity of the growth rate function on the space of primitive majors, and then prove it coincides with the value
given by Thurston's algorithm for rational majors. As a consequence, we get that the core entropy extends to a continuous function on ${\rm PM}(d)$,
establishing Theorem \ref{theorem:main}.

\begin{theorem}\label{theorem:continuity-simple}
The growth rate function $r : {\rm PM}(d)\to \mathbb{R}$ is continuous.
\end{theorem}
\begin{proof}
On the contrary, assume that there exists $\epsilon_0 > 0$ and a sequence $(m_N)$ of majors converging to $m$ such that $|r(m_N)-r(m)|>\epsilon_0$
for all $N$.
According to Proposition \ref{pro:limit-of-major}, there exists a subsequence which Hausdorff-converges to a critical portrait $\xi$, and $\xi$ induces $m$.
Moreover, by Lemma \ref{L:subseq}, by passing to a further subsequence (which we will still denote $(m_N)$ with abuse of notation) we can
assume that the
associated labeled wedges $\WWW_{N}$ converge to some labeled wedge $\WWW_\infty$.
By Lemma \ref{L:wplimit}, the limit wedge $\WWW_\infty$ is weakly periodic of type $\xi$. Let us denote as $\Gamma_N$ the infinite graph associated to $\mathcal{W}_N$, and $\Gamma_\infty$ the graph associated to $\mathcal{W}_\infty$.
As a consequence, we have for the growth rates
\[\text{$r(\G_{N})\to r(\G_\infty)$ as $N \to \infty$.}\]
Now, combining Lemma \ref{lem:limit-labeled wedge} and Proposition \ref{lem:common-growth-rate} we get
$$r(\Gamma_\infty) = r(\Gamma_\infty^{ND}) = r(\Gamma_\infty^Q)$$
and similarly, if $\Gamma_m$ denotes the infinite graph associated to the primitive major $m$ (note that $m$, being a primitive major, is trivially the limit of a constant family of primitive majors)
$$r(\Gamma_m) = r(\Gamma_m^{ND}) = r(\Gamma_m^Q)$$
Now, since $\Gamma_\infty$ is weakly periodic of type $\xi$, $\Gamma_m$ is weakly periodic of type $m$, and $m$ and $\xi$ are equivalent,
we have by Proposition \ref{pro:equi-portrait}
$$r(\Gamma_\infty^Q) = r(\Gamma_m^Q)$$
hence combining the previous equalities yields
$$r(\Gamma_\infty) = r(\Gamma_m)$$
which contradicts the assumption that $r(m_{N})\not\to r(m)$.
\end{proof}

To finish the proof of the main theorem, we need the following lemma.
\begin{lemma}\label{lem:entropy-continuity}
Let $\xi$ be a rational critical portrait of degree $d\geq2$. Then the logarithm of the growth rate $r(\xi)$ of the infinite graph $\G_\xi$ coincides with the core entropy of $\xi$:
\[h(\xi)=\log r(\xi).\]
\end{lemma}
\begin{proof}
Let $\xi$ be a rational critical portrait, $\Gamma_\xi$ its associated infinite graph, and $G_\xi := \G^Q_\xi$ the quotient graph of $\G_\xi^{ND}$.
By unraveling the definition, the matrix $A_\xi$ constructed in section \ref{sec:algorithm} is exactly the adjacency matrix of $G_\xi$. By Lemma \ref{lem:limit-labeled wedge} and Proposition \ref{lem:common-growth-rate}, the growth rate of $\G_\xi$ coincides with that of $G_\xi$. Moreover, by Lemma \ref{lem:adjacent-matrix}, the growth rate of $G_\xi$ coincides with the largest real eigenvalue of its adjacency matrix, that is the largest real eigenvalue of $A_\xi$. Thus, its logarithm is the core entropy $h(\xi)$.
\end{proof}

\begin{proof}[Proof of Theorem \ref{theorem:main}]
It follows directly from Lemma \ref{lem:entropy-continuity} and Theorem \ref{theorem:continuity-simple}.
\end{proof}

\section{Continuity of core entropy on the space of polynomials}

Let $d\geq2$ be an integer, and $f$ a complex polynomial of  degree $d$. The \emph{ filled-in Julia set} $K_f$ is the set of points
which do not escape to infinity under iteration, the \emph{Julia set} $J_f$ is the boundary of $K_f$ and the \emph{Fatou set} is $F_f:=\C\setminus J_f$.
 A point $c\in\C$ is called a \emph{critical point} of $f$ if $f'(c)=0$. The \emph{critical set} $ {\rm crit}(f)$  is defined to be $${\rm crit}(f)=\{c\in\C\mid f'(c)=0\},$$ and the \emph{postcritical set}  ${\rm post}(f)$ is defined to be
\begin{equation*}
  {\rm post}(f)=\ov{\{f^n(c) : {c\in{\rm crit}(f)}, {n\ge 1} \}}.
\end{equation*}

A polynomial is called \emph{\pf} if its postcritical set is finite.
Any \pf polynomial $f$ has a $f$-invariant tree $H_f$ containing the orbits of its critical points, called the \emph{Hubbard tree}, which captures the dynamics of the polynomial. Following Thurston, the \emph{core entropy} of $f$, denoted by $h(f)$, is defined to be the topological entropy of $f$ on its Hubbard tree, i.e.,\[h(f):=h_{top}(f,H_f).\]

In the previous part, we showed the continuity of the core entropy of rational critical portraits. As an application, we will prove the continuity of the core entropy of \pf polynomials of any given degree.

 Let $\PPP_d$ denote the parameter space of  monic centered polynomials of degree $d$. We say that a sequence of polynomials $(f_n)_{n\geq1}\subset \PPP_d$ \emph{converges} to $f\in\PPP_d$ if the coefficients of $f_n$ converge to the corresponding coefficients of $f$. The objective of this section is to prove the following result.

\begin{theorem}\label{core-entropy-polynomial}
Let $f_n,n\geq1$ and $f$ be \pf polynomials in $\PPP_d$. If $f_n\to f$ as $n\to\infty$, then $h(f_n)\to h(f)$.
\end{theorem}

We summarize the outline of the proof. Following Poirier \cite{Poi}, we associate to each polynomial $f_n$ (\emph{resp.} $f$) a rational formal critical portrait $\Theta_n=\{\FFF_n,\JJJ_n\}$ (\emph{resp.} $\Theta=\{\FFF,\JJJ\}$), called a (weak) \emph{critical marking} (see Section \ref{section:critical-marking} below). By Theorem \ref{entropy-algorithm-1}, we have
    \begin{equation}\label{eq:algorithm}
    h(\Theta_n)=h(f_n)\text{ for }n\geq1\ \text{ and }\ h(\Theta)=h(f).
    \end{equation}
Therefore, applying Theorem \ref{theorem:main}, one just needs to have a good choice of $\Theta_n$ and $\Theta$ such that $\Theta_n$ Hausdorff converge  to $\Theta$ as $n\to\infty$ and $\Theta$ is a weak critical marking for $f$. This is accomplished in Proposition \ref{Julia-type-converge} by studying continuity properties 
of external rays. 


\subsection{The dynamics of polynomials}
Let $f\in\PPP_d$. A point $z\in\C$ is called a \emph{preperiodic point} of $f$ if there exist integers $m\geq 0,n\geq1$ such that $f^{m}(z)=f^{m+n}(z)$. If $m=0$, the point $z$ is called \emph{periodic}. The minimal $m$ and $n$ with this property are called the \emph{preperiod} and \emph{period} of $z$ respectively.

Let $f$ be a polynomial in $\PPP_d$ with connected filled-in Julia set. By B\"{o}ttcher's Theorem, there exists a unique conformal isomorphism $\phi_f:\C\setminus \D\to \C\setminus K_f$ with $\phi_f$ tangent to the identity at $\infty$, such that the following diagram is commutative:
\begin{equation}\label{Bottcher}
\begin{array}{ccc}
\C\setminus K_f &\xrightarrow[]{\ \ f\ \ }  &\C\setminus K_f \\
\phi_f\Big\downarrow &&\Big\downarrow \phi_f  \\
\C\setminus \ov{\D} &  \xrightarrow[]{z\mapsto z^d} & \C\setminus \ov{\D}\vspace{-0.1cm}.\end{array}
 \end{equation}
The map $\phi_f$ is called the \emph{B\"{o}ttcher coordinate} of $f$. The \emph{external ray of argument $\theta$}, denoted by $R_f(\theta)$, is the image by $\phi_f^{-1}$ of the ray $\{z=re^{2\pi i\theta}\mid r>1\}$. We say that it \emph{lands} if the intersection
\[\bigcap_{r>1}\overline{\phi_f^{-1}((1,r]e^{2\pi i\theta})}\]
is a point, called the \emph{landing point} of $R_f(\theta)$. Since a power map sends radial lines to radial lines, the polynomial $f$ sends external rays to external rays. Set $U_f(\infty):=\C\setminus K_f$. The \emph{Green function} $G_f$ associated with $f$ is the harmonic function equal to $\log|\phi_f(z)|$ on $U_f(\infty)$ and vanishing on $K_f$. The number $s=G_f(z)\geq 0$ is called the \emph{potential} of $z\in\C$.

Now, we assume that $f$ is a \pf polynomial. Then the Fatou set of $f$ consists of attracting basins and all periodic points in $J_f$ are repelling. The filled-in Julia set $K_f$ is connected and locally-connected, and each bounded Fatou component is a Jordan domain. By B\"{o}ttcher theorem's, there is a system of Riemann mappings $$\Big\{\phi_U: \D\to U\,\Big|\, U \text{ bounded Fatou component}\Big\}$$ so that
each extends to a homeomorphism on the closure $\overline{\D}$, and the following diagram commutes for all $U$:
\begin{equation*}
 \begin{tikzpicture}
   \matrix[row sep=0.8cm,column sep=2.4cm] {
     \node (Gammai) {$ \overline \D $}; &
       \node (Gamma) {$\overline \D$}; \\
     \node (S2i) {$\overline U$}; &
       \node (S2) {$\overline{f(U)}$.}; \\
   };
   \draw[->] (Gamma) to node[auto=left,cdlabel] {\phi_{f(U)}} (S2);
   \draw[->] (S2i) to node[auto=right,cdlabel] {f} (S2);
   \draw[->] (Gammai) to node[auto=left,cdlabel] {\text{ power map }z^{d_{\tiny\mbox{$U$}}}} (Gamma);
   \draw[->] (Gammai) to node[auto=right,cdlabel] {\phi_U} (S2i);
 \end{tikzpicture}
 \end{equation*}
The image $\phi_U(0)$ is called the \emph{center} of the Fatou component $U$. It is easy to see that any center is mapped to a critical periodic point 
under finitely many iterations of $f$.
 The images in $U$ under $\phi_U$ of closed radial lines in $\ov{\D}$ are, by definition, the \emph{internal rays} of $U$.  As with external rays, the polynomial $f$ sends internal rays to internal rays.

Let $f$ be a \pf polynomial.
Then
any pair of points in the closure of a bounded Fatou component  can be joined in a unique way by a Jordan arc consisting of (at most two) segments of internal rays. We call such arcs \emph{regulated}.
Since $K_f$ is arc-connected, given two points $z_1, z_2\in K_f$, there is an arc $\gamma: [0,1]\to K_f$ such that $\gamma(0)=z_1$ and $\gamma(1)=z_2$. In general, we will not distinguish  between the map  $\g$ and its image. It is proved in \cite{DH} that such arcs can be chosen in a unique way so that  the intersection with the closure of a Fatou component is regulated. We still call such arcs regulated and denote them by $[z_1,z_2]$.
By \cite[Proposition 2.7]{DH}, the set
\[H_f:=\bigcup_{p,q\in {\rm post}(f)}[p,q] \]
 is a finite connected tree, called the \emph{Hubbard tree} of $f$.  A point $z\in J_f$ is called \emph{biaccessible} if there are at least two rays landing at $z$. The following result is well-known.

\begin{lemma}\label{biaccessible}
Let $f$ be a \pf polynomial. Then every biaccessible point in $J_f$ is mapped to the Hubbard tree of $f$ under finitely many iterates of $f$.
\end{lemma}

\begin{definition}[Core entropy of polynomials]\label{core-entropy}
The \emph{core entropy} of $f$, denoted by $h(f)$, is defined to be the topological entropy of the restriction of $f$ to its Hubbard tree $H_f$, i.e.,
\[h(f):=h_{top}(f,H_f).\]
\end{definition}

\subsection{Weak critical markings of \pf polynomials}\label{section:critical-marking}

In order to classify all \pf polynomials up to topological conjugacy, Poirier \cite{Poi} defined for any \pf polynomial a finite 
collection of combinatorial data, called a \emph{critical marking}, considering the set of rays landing at the critical points of $f$. 

In this section, we recall the definition of critical marking, and explain how it can be used to compute the core entropy. 
However, as we will see in section \ref{sec:not-closed}, the set of critical markings of \pf polynomials is not closed: indeed, if  a sequence of polynomials $\{f_n\}$ converges to a polynomial $f$ and the corresponding critical markings $\Theta_n$ of $f_n$ converge to $\Theta$, then $\Theta$ is not necessarily a critical marking of $f$. To solve this problem, we also introduce the more general notion of \emph{weak critical marking} (see also \cite{G}).

This construction requires the definition of supporting rays/arguments as follows.

\begin{definition}[supporting rays/arguments]\label{support-ray}
Let $U$ be a bounded Fatou component of a postcritically finite polynomial $f$, and let $z \in \partial U$ a point on its boundary.  The external rays landing at $z$ divide the
plane in finitely many regions. We label the arguments of these rays by $\theta_1,\ldots, \theta_k$ in counterclockwise cyclic order,
 so that $U$ belongs to the region delimited by $R(\theta_1)$
and $R(\theta_2)$ ($\theta_1=\theta_2$
if there is a single ray landing at $z$). The ray $R(\theta_1)$ (resp. $R(\theta_2)$) is called the \emph{left-supporting} (resp. \emph{right-supporting}) ray of $U$ at $z$, and the argument $\theta_1$ (resp. $\theta_2$) is called the \emph{left-supporting} (resp. \emph{right-supporting}) argument of $U$ at $z$.
\end{definition}

\subsubsection{Critical Fatou markings}\label{sec:critical-Fatou}
 Let $f$ be a \pf  polynomial of degree $d$, and let $U_1, \dots, U_n$ be its
critical Fatou components (i.e., the Fatou components containing a critical point).
Following Poirier \cite{Poi}, we now construct for each critical Fatou component $U$ a finite set $\Theta(U)$, whose elements are angles of external rays which land on the boundary of $U$.
Denote $\delta_U := \text{deg}(f|_U)$.
\begin{itemize}
\item \textbf{Case 1}: We first consider the case when $U$ is a periodic, critical Fatou component. Let
\[U\mapsto f(U)\mapsto\cdots\mapsto f^n(U)=U \]
be a critical Fatou cycle of period $n$. We will construct the associated set $\Theta(U')$ for every critical Fatou component $U'$ in this cycle simultaneously.
Let $z\in\partial U$ be a periodic point  with period less than or equal to $n$. Let $\theta$ denote the left-supporting argument of $U$ at $z$. Clearly, $\theta$ is periodic with period $n$. We call $\theta$ a \emph{preferred} angle for $U$.
 Note that this choice
naturally determines a left-supporting argument of each Fatou component $f^k(U)$ for $k\in\{0,\ldots,n-1\}$, which is called a \emph{preferred angle} of $f^k(U)$.
Let $U'$ be a critical Fatou component in the cycle and $\theta'$  its preferred angle. 
We now define $\Theta(U')$ as any set of $\delta_{U'}$ angles such that: 
\begin{itemize}
\item[(a)]
$\theta' \in \Theta'(U)$; 
\item[(b)]
the rays corresponding to the elements of $\Theta(U')$ land at $\delta_{U'}$ distinct points of $\partial U'$ and are inverse images
of $f(R(\theta'))$. 
\end{itemize}
\item \textbf{Case 2}: $U$ is a strictly preperiodic Fatou component. Let $k$ be the minimal number such that $U'=f^k(U)$ is a critical Fatou component. We may assume that $\Theta(U')$ is already chosen, according to the previous case. Choose an angle $\theta' \in \Theta(U')$. We define $\Theta(U)$ to be the set of arguments of the $\delta_{U}$  rays landing at $\delta_U$ distinct points of $\partial U$ that are $k$-th inverse images
of $R(\theta')$.
\end{itemize}

Let $f$ be a \pf polynomial. Then a \emph{weak critical Fatou marking} is a collection
$$\mathcal{F} = \{ \Theta(U_1), \dots, \Theta(U_s) \}$$
as given by the above construction such that the convex hulls in $\ov{\D}$ of $\Theta(U_1),\ldots,\Theta(U_n)$ have pairwise disjoint interiors, where $U_1, \dots, U_n$ are the critical Fatou components of $f$. Weak critical Fatou markings are not uniquely determined by $f$, and there are finitely many choices.
If all angles which appear in $\mathcal{F}$ are left-supporting ones, then we call $\mathcal{\FFF}$ a \emph{critical Fatou marking}, which is the original object considered by Poirier.

As an example, we consider the cubic polynomial $f(z)=z^3+\frac{3}{2}z^2$.  The critical point $z = 0$ is fixed, and the other critical point $z = -1$ is mapped to a repelling fixed point $z = 1/2$ (see Figure \ref{intersect}).
\begin{figure}[http]
\begin{tikzpicture}
\node at (0,0) {\includegraphics[width=6cm]{right.pdf}};
\node at (-1.75,2){\footnotesize{$\RRR(\frac{1}{3})$}};
\node at (-1.75,-2){\footnotesize{$\RRR(\frac{2}{3})$}};
\node at (3.45,0){\footnotesize{$\RRR(0)$}};
\node at(1,0){$U$};
\end{tikzpicture}
 \caption{The Julia set of $f(z)=z^3+\frac{3}{2}z^2$}\label{intersect}
\end{figure}
Then there is only one critical Fatou component $U$, which contains $0$, and the point $z$ of least period on $\partial U$ is $z = 1/2$, hence the preferred angle is $\theta' = 0$. Then we have two choices for a weak critical Fatou marking of $f$: namely, $\mathcal{F}=\{\ \Theta(U)=\{0,1/3\}\ \}$ is a weak critical Fatou marking, but not a critical Fatou marking, since $R(1/3)$ is not left-supporting for $U$, while $\mathcal{F}=\{\ \Theta(U)=\{0,2/3\}\ \}$ is a critical Fatou marking of $f$.

\subsubsection{Critical Julia markings}\label{sec:critical-julia-marking} Let $c$ be a critical point which lies in the Julia set of $f$. Then a \emph{critical Julia leaf} landing at $c$ is a finite subset $\Theta$ of the circle of cardinality $\geq 2$ such that:
\begin{enumerate}
\item
for each $\theta \in \Theta$, the external ray with angle $\theta$ lands at $c$;
\item
all rays $R(\theta)$ with $\theta \in \Theta$ are mapped by $f$ to the same ray.
\end{enumerate}

A \emph{weak critical Julia marking} of $f$ is a finite collection
$$\mathcal{J} = \{ \Theta_1(c_1), \dots, \Theta_m(c_m) \}$$
where:

\begin{enumerate}
\item each $\Theta_i(c_i)$ is a critical Julia leaf landing at $c_i$;
\item the set $\{c_1, \dots, c_m\}$ equals the set of all critical points of $f$ which lie in the Julia set (however, the $c_i$ need not be distinct!)
\item any two of the convex hulls 
in the closed unit disk of $\Theta_1(c_1),\ldots,\Theta_m(c_m)$  either are disjoint or intersect at one point on $\partial\D$;
\item for each critical point $c\in J_f$, we have the formula
\[{\rm deg}(f|_c)-1=\sum_{c_j=c}\big(\#\Theta_j(c_j)-1\big). \]
\end{enumerate}

Once again, weak critical Julia markings are not uniquely determined by $f$, and there are finitely many choices.
If all $c_i$ are distinct, then we call $\mathcal{J}$ a \emph{critical Julia marking}. Critical Julia markings are the original combinatorial objects
defined by Poirier, while we relax the definition by allowing the same critical point in the Julia set to appear with multiplicity.

\begin{figure}
\begin{tikzpicture}
\node at (0,0) {\includegraphics[width=6.5cm]{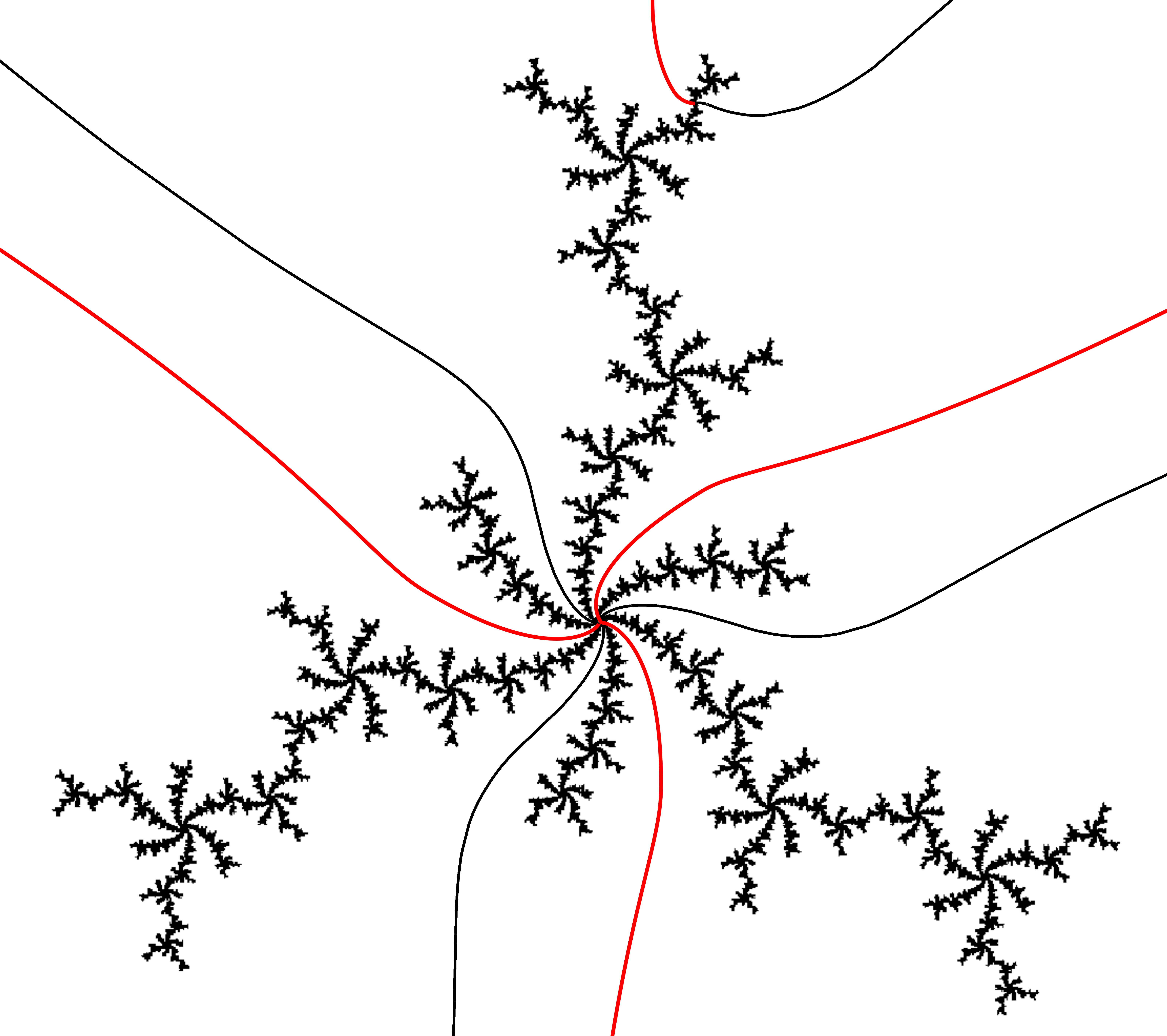}};
\node at (2.25,3){\footnotesize{$\frac{11}{72}$}};
\node at (0.25,3){\footnotesize{$\frac{17}{72}$}};
\node at (3.5,1.25){\footnotesize{$\frac{17}{216}$}};
\node at (3.5,0.25){\footnotesize{$\frac{11}{216}$}};
\node at (-3.5,2.5){\footnotesize{$\frac{83}{216}$}};
\node at(-3.5,1.5){\footnotesize{$\frac{89}{216}$}};
\node at(-1,-2.5){\footnotesize{$\frac{155}{216}$}};
\node at(0.5,-2.5){\footnotesize{$\frac{161}{216}$}};
\end{tikzpicture}
 \caption{The Julia set of $f_c(z)=z\mapsto z^3+0.22036+1.18612 i$.}\label{portrait}
\end{figure}

To show the non-uniqueness, let us consider the following example, which comes from \cite{G}. We consider the \pf polynomial $f_c(z)=z^3+c$ with $c\approx 0.22036+1.18612 i$. The critical value $c$ receives two rays with arguments $11/72$ and $17/72$. Then, $$\Theta:=\large\{\ \Theta_1(0):=\left\{11/216,83/216\right\},\Theta_2(0):=\left\{89/216,161/216\right\}\ \large\}$$
is a weak critical marking, but not a critical marking, of $f_c$, and
\[\Theta:=\large\{\ \Theta(0):=\{11/216,83/216,155/216\}\ \large\}\]
is a critical marking of $f_c$ (see Figure \ref{portrait}).

\begin{definition}\label{weak-critical-marking}
A weak critical marking of $f$ is a collection
\begin{equation}\label{eq:weak-critical-marking}
\Theta = \{ \mathcal{F}, \mathcal{J} \}
\end{equation}
where $\mathcal{F}$ is a weak critical Fatou marking of $f$ and $\mathcal{J}$ is a weak critical Julia marking of $f$, such that the convex hulls in $\ov{\D}$ of the elements of $\Theta$ have pairwise disjoint interiors.
\end{definition}

Note that a weak critical marking $\Theta$ of a \pf polynomial is 
a rational critical portrait (but not necessarily a primitive major! See \cite{Poi}, Example 2.7).
Thus, we will denote as $h(\Theta)$ the core entropy associated to this critical portrait.
Any such polynomial admits at least one, and in general finitely many weak critical markings.
If $\mathcal{F}$ and $\mathcal{J}$ are actually a critical Fatou marking and a critical Julia marking (as opposed to weak ones), then we call $\Theta$ a \emph{critical marking} of $f$.

\subsection{Convergence of external rays}

In the proof of Theorem \ref{core-entropy-polynomial}, a result about the convergence of external rays (Lemma \ref{ray-convergence-2}) will play an important role. The aim of this section is to prove this convergence result based on a sequence of lemmas (some of which are well-known).

\begin{lemma}\label{continue}
Let $f\in \PPP_d$ and $z$ be a repelling preperiodic point of $f$ such that the forward orbit of $z$ avoids the critical points of $f$. Then there exists a neighborhood $\Lambda$ of $f$ in $\PPP_d$ and a holomorphic map $\xi_z:\Lambda\to \C$ such that $\xi_z(f)=z$ and $\xi_z(f')$ is the unique repelling preperiodic point of $f'$ near $z$ with the same preperiod and period as $z$ for all $f'\in \Lambda$. The point $\xi_z(f')$  is called the \emph{continuation of $z$ at $f'$}.
\end{lemma}

The proof follows directly from the implicit function theorem. Let now $\{S_n\}\subset \C$ be a sequence of sets. We denote as
\[\limsup S_n\]
the set of points $z\in\C$ such that every neighborhood of $z$ intersects infinitely many $S_n$. It follows immediately from the definition that $\limsup S_n$ is closed.

\begin{lemma}[Goldberg-Milnor \cite{GM}]\label{perterb-ray}
Consider a polynomial $f\in \PPP_d$ and an external ray $R_f(\theta)$ which lands at a repelling preperiodic point $z$ such that the orbit of $z$ avoids the critical points of $f$. Then $R_{f'}(\theta)$ lands at the continuation of $z$ at $f'$, for all $f'$ in a sufficiently small neighborhood of $f$. Moreover, if $f_n\to f$ as $n\to\infty$, then $\limsup_{n} \ov{R_{f_n}(\theta)}=\overline{R_f(\theta)}$.
\end{lemma}

Assume that $f_n,n\geq 0$ are polynomials in $\PPP_d$ with connected Julia set.
 For each $n\geq0$, we simply denote $J_{f_n},K_{f_n}$ by $J_n,K_n$ respectively, the external ray $R_{f_n}(\theta)$ by $R_n(\theta)$ for all $\theta\in\R/\Z$, the infinite Fatou component $U_{f_n}(\infty)$ by $U_n(\infty)$, and the B\"{o}ttcher coordinate $\phi_{f_n}$
     given in \eqref{Bottcher} by $\phi_n$. The following result is well-known.

\begin{lemma}\label{Riemann-mapping-convergent}
Let $f_n,n\geq 0$ be  polynomials in $\PPP_d$ with connected Julia set such that $f_n\to f_0$ as $n\to\infty$. Then the inverse $\psi_n$ of the B\"{o}ttcher coordinate $\phi_n$ converges to $\psi_0:=\phi_0^{-1}$ uniformly on any compact subset of $\C\setminus\ov{\D}$.
\end{lemma}

\begin{proof} Let $\psi_n : \mathbb{C} \setminus K(f_n) \to \mathbb{C} \setminus \overline{\mathbb{D}}$ be the Riemann map of $K(f_n)$. Since all $K(f_n)$ are uniformly bounded, then all the images of $\psi_n$ contain a ball around $\infty$ of uniform radius. Hence, the family $(\psi_n)$ is precompact: let $\psi_0$ be any limit. Each $\psi_n$ satisfies the B\"ottcher equation
$\psi_n(z^d) = f_n(\psi_n)$, hence by taking the limit one gets $\psi_0(z^d) = f_0(\psi_0)$, so $\psi_0$ is the B\"ottcher map for $f_0$.
\end{proof}

\begin{lemma}\label{limit-of-ray}
Let $f_n,n\geq 0$ be polynomials in $\PPP_d$ with connected Julia set such that $f_n\to f_0$ as $n\to\infty$. For each argument $t$ and any sequence of arguments $\sigma=\{t_n\}_{n\geq1}$ converging to $t$, let us denote $B_\sigma(t):=\limsup\ov{R_n(t_n)}$ and $K_\sigma(t):=B_\sigma(t)\cap K_0$, and also denote $\tau(\sigma):=\{\tau(t_n)\}_{n\geq1}$ (it converges to $\tau(t)$). Then we have
\begin{enumerate}
\item  the intersection of $B_\sigma(t)$ and $U_0(\infty)$ is $R_0(t)$, so that $B_\sigma(t)=R_0(t)\cup K_\sigma(t)$;
\item  the sets $B_\sigma(t)$ and $K_\sigma(t)$ are connected, $f_0(B_\sigma(t))\subset B_{\tau(\sigma)}(\tau(t))$ and $f_0(K_\sigma(t))\subset K_{\tau(\sigma)}(\tau(t))$.
\end{enumerate}
\end{lemma}
\begin{proof}
 (1) On the one hand, let $z_n\in\ov{R_n(t_n)},n\geq1,$ converge to $z\in B_\sigma(t)$, and the potential of $z_n$ be $s_n$. By choosing a subsequence if necessary, we assume that $s_n\to s\geq 0$ as $n\to \infty$. It is known that the Green functions $G_n(z)$ uniformly converge to $G_0(z)$ on $\C$ (\cite[Proposition 8.1]{DH}), so $z\in U_0(\infty)$ if and only if $s>0$. In the case of $s>0$, by Lemma \ref{Riemann-mapping-convergent}, the points $z_n=\psi_n(e^{s_n+2\pi it_n})$ converge to $z=\psi_0(e^{s+2\pi it})\in R_0(t)$. On the other hand, given any $s>0$, by Lemma \ref{Riemann-mapping-convergent}, we have
    \[R_n(t_n)\ni\psi_n(e^{s+2\pi i t_n})\to \psi_0(e^{s+2\pi it})\in R_0(t)\]
Since $s$ is arbitrary, it follows that $\ov{R_0(t)}$ is contained in $B_\sigma(t)$.

  \noindent (2) Let $x$ be a point of $\ov{R_0(t)}\cap K_0$, which belongs to $K_\sigma(t)\subset B_\sigma(t)$ by (1). Let now $y$ be another point of $B_\sigma(t) \cap K_0$: by definition, there exists a sequence $y_{n_k}$
such that $y_{n_k} \in R_{n_k}(t_{n_k})$ and $y_{n_k} \to y$. By applying (1) to this subsequence, there exists a further subsequence (which we still denote  by $x_{n_k}$) of points which converge to $x$
and such that $x_{n_k} \in R_{n_k}(t_{n_k})$. Let us denote as $c_{n_k}$ the segment of the ray $R_{n_k}(t_{n_k})$ connecting $x_{n_k}$ and $y_{n_k}$, and let $c$ be a Hausdorff limit of the segments $c_{n_k}$.
Then by construction the set $c$ is a connected, compact set which contains $x$ and $y$, and it is also a subset of $B_\sigma(t)$, proving that $B_\sigma(t)$ is connected.
Note that $x$ and $y$ belong to $K_0$, so the potentials of $x_{n_k}$ and $y_{n_k}$ with respect to $f_{n_k}$ converge to $0$. This implies that the limit $c$ of $c_n$ belongs to $B_\sigma(t)\cap K_0=K_\sigma(t)$, proving that $K_\sigma(t)$ is connected.

Since $f_0(R_0(t))=R_0(\tau(t))$, it remains to show that $f_0(K_\sigma(t))\subset K_{\tau(\sigma)}(\tau(t))$. Let $z\in K_\sigma(t)$. Then there exist $z_n\in R_n(t_n)$ with potential $s_n$ such that $z_n\to z$ and $s_n\to 0$ as $n\to \infty$. Since $f_n$ uniformly converge to $f_0$, then $w_n:=f_n(z_n)$ converge to $w:=f_0(z)$ as $n\to\infty$. On the one hand, note that $w_n\in R_n(\tau(t_n))$, so $w\in B_{\tau(\sigma)}(\tau(t))$. On the other hand, the potentials of $w_n$ are $ds_n$, converging to $0$, so $w\in K_0$. It follows that $w=f(z)\in K_{\tau(\sigma)}(\tau(t))$.
\end{proof}

The next lemma comes directly from \cite[Lemma 6.3]{D}.

\begin{lemma}\label{attracting-basin-converge}
Let $f_n,n\geq0,$ be \pf polynomials in $\PPP_d$ such that $f_n\to f_0$ as $n\to\infty$. Let $U_0$ be a Fatou component of $f_0$. Then the center of $U_0$ is contained in a Fatou component of $f_n$, denoted by $U_n$, for all sufficiently large $n$. Furthermore, any given compact subset of $U_0$ is contained in $U_n$ for all sufficiently large $n$. The Fatou component $U_n$  is called the \emph{deformation of $U_0$ at $f_n$}.
\end{lemma}

\begin{lemma}\label{center-converge}
Let $f_n,n\geq 0,$ be \pf polynomials in $\PPP_d$ such that $f_n\to f_0$ as $n\to\infty$. Let $U_0$ be a Fatou component of $f_0$, and $U_n$ the deformation of $U_0$ at $f_n$ for each sufficiently large $n$. Then the centers of $U_n$ converge to that of $U_0$ as $n\to\infty$, and ${\rm deg}(f_n|_{U_n})={\rm deg}(f_0|_{U_0})$ for all sufficiently large $n$. Furthermore, for any preperiodic point $z\in\partial U_0$,  there is a unique point $z_n\in\partial U_n$, having the same preperiod and period as $z$, such that $z_n\to z$ as $n\to\infty$. The point $z_n$  is called the \emph{continuation of $z$ at $\partial U_n$}.
\end{lemma}
\begin{proof}
Let $x_n$ be the center of $U_n$.
If $x_0$ is periodic, the continuation $y_n$ of $x$ at $f_n$ is an attracting periodic point contained in $U_n$ (by the Implicit Function Theorem and Lemma \ref{attracting-basin-converge}). Hence $z_n=p_n$.
Let us now deal with the preperiodic case by induction. Let us assume that $f_n(x_n)\to f_0(x_0)$ as $n\to\infty$: we need to show that
$x_n\to x_0$ and ${\rm deg}(f_n|_{U_n})=\text{deg}(f_0|_{U_0})$ as $n\to\infty$.
Set $\delta:={\rm deg}(f_0|_{U_0})$. By Rouch\'{e}'s theorem, any given small neighborhood of $p_0$  contains exactly $\delta$ preimages by $f_n$ of $f_n(x_n)$ (counting with multiplicity) for every sufficiently large $n$. Note that all these preimages belong to $U_n$ by Lemma \ref{attracting-basin-converge}, and are the centers of some Fatou component of $f_n$. So these preimages must coincide with $x_n$. It follows that $x_n\to x_0$ as $n\to\infty$ and $\text{deg}(f_n|_{U_n})=\delta$ for all sufficiently large $n$.

For the remaining result of this lemma, we first assume that $z$ is periodic. Then $z$ is repelling because $f_0$ is postcritically finite. In this case, the conclusion holds by Goldberg and Milnor's proof in \cite[Appendix B]{GM}.
Now, let $z\in\partial U_0$ be a preperiodic point. Set $v:=f_0(z)\in \partial f(U_0)$. Inductively, we assume that $v_n$ is the unique preperiodic point of $f_n$ in $\partial f_n(U_n)$ such that $v_n$ has the same preperiod and period as $v$, and $v_n\to v$ as $n\to\infty$.  Since $f_n$ uniformly converges to $f_0$, given any small disk neighborhood $W_z$ of $z$, there is a disk neighborhood $V_v$ of $v$ such that the component of $f^{-1}_n(V_v)$ that contains $z$, denoted by $D_{n,z}$, belongs to $W_z$, for all sufficiently large $n$ and $n=0$. Given any sufficiently large $n$, choose a point $a_n\in D_{n,z}\cap U_n$ and set $b_n:=f_n(a_n)$. Then $b_n\in V_v\cap f_n(U_n)$.
By the inductive assumption, the point $v_n$ belongs to $\partial f_n(U_n)\cap V_v$. One can then choose an arc $\g_n\subset f_n(U_n)\cap V_v$ joining $b_n$ and $v_n$. Lifting $\g_n$ by $f_n$ with the starting point $a_n$, we get an arc $\wt{\g}_n\subset D_{n,z}\cap U_n$. Its ending point, denoted by $z_n$,  belongs to $\partial U_n$ and satisfies that $f_n(z_n)=v_n$.
By the argument above, we in fact proved that for any point $z'\in\partial U_0$ with $f_0(z')=v$, and any small neighborhood $W_{z'}$ of $z'$, there exists a point $z_n'\in\partial U_n$ with the property that $z_n'\in W_{z'}$ and $f_n(z_n')=v_n$ for all sufficiently large $n$. Since $\text{deg}(f_n|_{U_n})=\text{deg}(f_0|_{U_0})$, the points which have the same properties as $z_n'$ are unique.
This completes the proof of the lemma.
\end{proof}

\begin{lemma}\label{lem:internal-ray}
Let $f_n,n\geq 0,$ be \pf polynomials in $\PPP_d$ such that $f_n\to f_0$ as $n\to\infty$. Let $U_0$ be a Fatou component of $f_0$, and $U_n$ the deformation of $U_0$ at $f_n$ for each large $n$. Suppose $I_n$ is a preperiodic internal ray of $f_n$ in $U_n$ with fixed preperiod $k\geq0$ and period $p\geq1$. If the landing point $z_n$ of $I_n$ converges to $z$, then $\limsup_{n\to\infty}I_n=I,$ where $I$ is the internal ray of $f$ in $U_0$ landing at $z$.
\end{lemma}

\begin{proof}
If $I$ is periodic, the conclusion holds by Goldberg and Milnor's proof in \cite[Appendix B]{GM}. By induction on $k$, it then suffices to prove $\limsup I_n=I$ provided that $\limsup f_n(I_n)=f(I)$. Since $f_n\to f$, we can choose  B\"{o}ttcher coordinates $\varphi_0$ of $U_0$ and $\varphi_n$ of $U_n$ such  that $\varphi_n^{-1}:D\to U_n$ converge
uniformly on compact sets to $\varphi^{-1}_0:D\to U_0$. It follows that $I':=\limsup I_n\cap U_0$ is an internal ray of $U$. On the other hand, note that $\limsup I_n\cap \partial U_0$ is compact, connected and contains the point $z$. The map $f$ sends $\limsup I_n\cap \partial U_0$ into the set $\limsup f_n(I_n)\cap\partial f(U_0)$, which is by induction a singleton. Then we get $\limsup I_n\cap \partial U_0=\{z\}$, and hence $I'=I$.
\end{proof}

\begin{lemma}\label{key}
Let $f_n,n\geq 0,$  be \pf polynomials in $\PPP_d$ such that $f_n\to f_0$ as $n\to\infty$. Let $U_0$ be a Fatou component of $f_0$, and $U_n$ the deformation of $U_0$ at $f_n$.  If $\theta$ is the left-supporting (resp. right-supporting) angle of $U_0$ at a periodic point $z$, then $\theta$ is also the left-supporting (resp. right-supporting) angle of $U_n$ at $z_n$ for all  large $n$, where $z_n$ denotes the continuation of $z$ at $f_n$.
\end{lemma}
\begin{proof}
We just prove this lemma in the case that $\theta$ is a left-supporting angle for $U_0$. The proof of the right-supporting case is exactly the same.
 Let $\theta_1,\ldots,\theta_s$ be the external angles associated with $z$ in the counterclockwise direction with $\theta_1=\theta$. In this case, all $\theta_1,\ldots,\theta_s$ are periodic with a common period, and $z$ is a repelling periodic point. By Lemma \ref{center-converge}, the continuation $z_n$ of $z$ at $f_n$ belongs to $\partial U_n$ for all large $n$, and Lemma \ref{perterb-ray} implies that the external rays of $f_n$ with arguments $\theta_1\ldots,\theta_s$ land at $z_n$.

Pick a point $p\in U_0$. We denote by $W$ the component of $\C\setminus( R_{0}(\theta_1)\cup R_{0}(\theta_2))$ that contains $p$.
Since $\limsup R_n(\theta_i)=R_0(\theta_i)$
for all $i=1,\ldots,s$, then, for each sufficiently large $n$, there exists a unique component of $\C\setminus( R_n(\theta_1)\cup R_n(\theta_2))$ that contains $p$, which we denote by $W_n$. Note that $p\in U_n$ and $U_n$ is contained in a component of $\C\setminus( R_n(\theta_1)\cup R_n(\theta_2))$, so $U_n\subset W_n$ for all sufficiently large $n$. We denote by $(\theta_1,\theta_2)$ the set of arguments we meet when traveling on $\R/\Z$ from $\theta_1$ to $\theta_2$ in the counterclockwise direction.

By contradiction, and passing to a subsequence if necessary, one can assume that $\theta_1$ is not the left-supporting angle of $U_{n}$ at $z_n$ for all sufficiently large $n$. For each $n$, we denote by $\eta_n$ the left-supporting angle of $U_{n}$ at $z_{n}$. By the argument in the last paragraph, each $\eta_n$ belongs to $(\theta_1,\theta_2)$. Note also that each $\eta_n$ has the same period as $\theta_1$ so, by choosing a subsequence if necessary, one can assume $\eta_n=\eta\in (\theta_1,\theta_s)$ for all sufficiently large $n$. But then, by Lemma \ref{perterb-ray}, the ray $ R_0(\eta)$ also lands at $z$, contradicting the fact that $\theta_1$ is the left-supporting angle for $U_0$ at $z$.
\end{proof}

\begin{lemma}\label{density}
Let $f$ be a \pf polynomial, and $S\subset J_f$ be a connected compact set with more than one point. Let $[z,w]$ denote the regulated arc in $K_f$ joining $z\not=w\in S$.
\begin{enumerate}
\item Every component of $S\cap [z,w]$ is an arc or a point in $J_f$; and every component of $[z,w]\setminus S$ is the union of two internal rays of a Fatou component $U$.
\item  If $(a,b)$ is a component of $[z,w]\setminus S$ with $b\not\in\{z,w\}$, then $b$ is a preperiodic point in the boundary of the Fatou component containing $(a,b)$.
\item If $[z,w]\subset S$, then the open arc $(z,w)$ contains a preperiodic point whose forward orbit avoids the critical points of $f$.
\end{enumerate}
\end{lemma}
\begin{proof}
(1) The first conclusion is obvious because $[z,w]$ is an arc. To prove the second one, let $(a,b)$ be a component of $[z,w]\setminus S$. Then there exists a bounded component $D$ of $\C\setminus([z,w]\cup S)$ such that $\partial D$ contains $(a,b)$. Note that $D$ belongs to the interior of $K_f$, so it belongs to a Fatou component $U$. It follows that
    $\partial D\setminus (a,b)\subset \ov{U}\cap S\subset \partial U$. Hence, $a,b\in \partial U$ and $(a,b)$ is the union of the two internal rays in $U$ landing at $a$ and $b$.

\noindent (2) In this case, $b$ is a biaccessible point, i.e., there are at least two external rays landing at $b$, and is in the boundary of a Fatou component according to (1). By Lemma \ref{biaccessible}, all its sufficiently high iterates by $f$ are intersections of periodic Fatou components and the Hubbard tree. Since there are only finitely many such points, then $b$ is preperiodic.

\noindent(3)
As $f$ is \pf,
then it is expanding in a neighborhood of $J_f$ in the sense that,
given a neighborhood $W$ of $J_f$, there exist constants $\lambda>1$ such that for any arc $\g\subset J_f$ with $f^n:\g\to\C$  injective,
\begin{equation}\label{eq:expanding}
{\rm length}(f^n(\g))\geq \lambda^n {\rm length}(\g),
\end{equation}
where ${\rm length}(\cdot)$ denotes the length of arcs in the canonical orbifold metric of $f$
(see \cite[Section 4]{DH}, \cite[Section 19]{Mil} and \cite[Section A.3]{Mc}).

We denote by $\AAA$ the set of open regulated arcs $(c,\xi)$ satisfying the conditions
\begin{enumerate}
\item $c$ is a critical point of $f$ and $[c,\xi]\subset H_{f}\cap J_f$;
\item $(c,\xi)$ avoids the postcritical points of $f$ and the branching points of $H_{f}$;
\item ${\rm length((c,\xi))=\kappa}$, where $\kappa$ is a sufficiently small universal constant.
\end{enumerate}
It is clear that $\AAA$ contains finitely many elements.

We claim that the preperiodic points whose forward orbit avoids the critical points of $f$ are dense in each member of $\AAA$.
Given an arc $\gamma$, a point $z$ on $\gamma$, and $\delta > 0$, we say that $z$ is $\delta$-\emph{contained} in $\gamma$
if $\gamma$ contains an open arc of length $2 \delta$ with center $z$. To represent such an arc, we use the notation
\[D_{\delta}^{\g}(z):=\{w\in\g \mid {\rm length}([z,w])< \delta \}.\]
To prove the claim, let $\g_1$ be any element of $\AAA$, and pick a point $ a\in \g_1$ and $\epsilon>0$.
Let us now choose a number $\delta_1<\epsilon/2$ such that $a$ is $2\delta_1$-contained in $\g_1$.  Since $f$ is expanding, the forward iterates of any open segment in $\g_1$ will eventually contain a critical point of $f$. It follows that there exists a sufficiently large integer $n_1$ with $\kappa/\lambda^{n_1}<\delta_1$ and a segment $[z_1,w_1]\subset D_{\delta_1}^{\g_1}(a)$ such that $[z_2,w_2]=f^{n_1}([z_1,w_1])$ belongs to an element of $\AAA$, denoted $\g_2$. Let $\delta_2>0$ such that
$z$ is $\delta_2$-contained in $\gamma_2$ 
for every $z\in [z_2,w_2]$. By shrinking $[z_1, w_1]$ if necessary, one can find an integer $n_2$ with $\kappa/\lambda^{n_2}<\delta_2$, such that $[z_3,w_3]:=f^{n_2}([z_2,w_2])$ is contained in an element of $\AAA$, denoted $\g_3$.
Repeating this process $N:=\#\AAA$ times, we obtain the segments $[z_i,w_i]$ and the elements $\g_i$ of $\AAA$ for $i=1,\ldots,N+1$, and the numbers $n_i, \delta_i$ for $i=1,\ldots, N$, such that
\begin{itemize}
\item $[z_i,w_i]\subset\g_i\in\AAA$;
\item every $z\in[z_i,w_i]$ is $\delta_i$-contained in $\g_i$;
\item $\kappa/\lambda^{n_i}<\delta_i$, and
\item $f^{n_i}[z_i,w_i]=[z_{i+1},w_{i+1}]$.
\end{itemize}
For each $i\in\{1,\ldots,N\}$, we denote by $\beta_i$ the lift of $\g_{i+1}$ by $f^{n_i}$ that contains $[z_i,w_i]$. We claim that $\beta_i\subset D_{\delta_i}^{\g_i}(z_i)\subset \g_i$.
Since $f$ is uniformly expanding on $J_f$ and by the choice of $n_i,\delta_i$, the length of $\beta_i$ satisfies
 \[{\rm length}(\beta_i)\leq {\rm length}(\g_{i+1})/\lambda^{n_i}=\kappa/\lambda^{n_i}<\delta_i.\]
So it is enough to prove that $\beta_i\subset \g_i$. On the contrary, there must be a point $p\in \beta_i\cap \g_i$ which is a branch point of $\g_i\cup \beta_i$. By property (2) in the construction of $\AAA$, the first $n_i$ terms in the orbit of  $p$ contain no critical points of $f$.
Then $f^{n_i}(p)$ is a branch point of $f^{n_i}(\beta_i) \cup f^{n_i}(\gamma_i)$; now, $f^{n_i}(\beta_i) = \gamma_{i+1}$ is a subset of $H_f$, and moreover $\gamma_i \subseteq H_f$
so also $f^{n_i} (\gamma_i) \subseteq H_f$.
Thus, $f^{n_i}(p)$ is a branch point of the Hubbard tree $H_f$,
which contradicts property (2) and completes the proof of the claim.
Since $\#\AAA=N$, there exist $i<j\in\{1,\ldots,N+1\}$ such that $\g_i=\g_j$. Denote by $\g_i'$ the pullback of  $\g_j=\g_i$ along the orbit from $[z_i,w_i]$ to $[z_j,w_j]$. It follows from the claim above that $\g_i'\subset D_{\delta_i}^{\g_i}(z_i)\subset \g_i$. Then the attracting map $$(f^{n_i})^{-1}:\g_i\to\g_i'\subset \g_i$$ has a fixed point. Hence $D_{\delta_i}^{\g_i}(z_i)$ contains a periodic point, which is disjoint from the orbits of the critical points of $f$ by $(2)$. Consequently, $D_{\delta_1}^{\g_1}(z_1)\subset D^{\g_1}_\epsilon(a)$ contains a preperiodic point whose orbit avoids the critical points of $f$. Note that $\g_1,a\in\g_1$ and $\epsilon$ are all arbitrary, so the claim is proven.

 Since $[z,w]\subset J_f$, by shrinking $[z,w]$ if necessary,  each of $z,w$ receives at least two rays of $f$. By Lemma \ref{biaccessible}, $z$ and $w$ are eventually  mapped into the Hubbard tree by iterations of $f$. By shrinking $[z,w]$ again if necessary, one can assume that $[z',w']:=f^n([z,w]) \subset H_f\cap J_f$.
 Since $f$ is expanding, some iteration of $[z',w']$ must contain a critical point of $f$, and hence intersect some element of $\AAA$. It follows from the claim above that $[z,w]$ contains preperiodic points whose forward orbits avoid the critical points of $f$.
\end{proof}

Using the previous lemmas, we can prove the following convergence result.

\begin{lemma}\label{ray-convergence-2}
Let $f_n,n\geq0,$ be \pf polynomials in $\PPP_d$ such that $f_n\to f_0$ as $n\to\infty$. If the angles $\theta_n$ converge to an angle $\theta$, then $\limsup \ov{R_n(\theta_n)}=\ov{R_0(\theta)}$, and the landing points of $R_n(\theta_n)$ converge to that of $R_0(\theta)$.
\end{lemma}
\begin{proof}
Note that if the first conclusion holds, then the second one follows directly. So we just need to prove  $\limsup \ov{R_n(\theta_n)}=\ov{R_0(\theta)}$.
We follow the notation of Lemma \ref{limit-of-ray}.  Set $\sigma:=\{\theta_n\}_{n\geq1}$. It is enough to show that $K_\sigma(\theta)$ is a singleton.

If this is not the case, then, by Lemma \ref{limit-of-ray} the set $K_\sigma(\theta)$ is connected and contains a point $w$ distinct from the landing point $z$ of $R_0(\theta)$ which belongs to $K_\sigma(\theta)$. Moreover,
$K_\sigma(\theta)$ is contained in $J_0$: indeed, if there exists $z \in K_\sigma(\theta)$ which belongs to a Fatou component $U$, then by Lemma \ref{attracting-basin-converge}
it also belongs to its deformation $U_n$ for $n$ large, hence it cannot be an accumulation point of the rays $R_n(\theta_n)$.
Let $[z,w]$ denote the regulated arc in $K_0$.

In the case of $[z,w]\not\subset K_\sigma(\theta)$, the segment $[z,w]$ passes through a Fatou component $U$ of $P$. We choose an arc $\G$ separates $z,w$ as follows.
Pick $x,x'$ in different component of $\partial U\setminus [z,w]$ such that their orbits avoid the critical points of $P$. We denote $R_0(t),R_0(t')$ the external rays landing at $x,x'$ respectively, and $I,I'$ the internal rays in $U$ landing at $x,x'$ respectively. The arc $\G$ is defined as
\[\G:=R_0(t)\cup \ov{I}\cup\ov{I'}\cup R_0(t').\]
It clearly separates $z$ and $w$.
By Lemma \ref{perterb-ray}, the rays $R_n(t),R_{n}(t')$ land at the continuation $x_n,x_n'$ of $x,x'$ respectively, and they belong to the boundary of $U_n$ by Lemma \ref{center-converge}, where $U_n$ is the deformation of $U$ at $f_n$. We thus obtain an arc
$\G_n:=R_{n}(t)\cup \ov{I_n}\cup\ov{I'_n}\cup R_{n}(t')$ for each sufficiently large $n$,
with $I_n,I_n'$ the internal rays in $U_n$ landing at $x_n,x_n'$ respectively.
In the case of $[z,w]\subset K_\sigma(\theta)$ ($\subset J_f$), by Proposition \ref{density}.(3), there exists a preperiodic point $b\in(z,w)$ such that its forward orbit avoids the critical points of $f_0$. We pick two rays $R_0(t),R_0(t')$ landing at $b$ such that the simple curve $\G:=R_0(t)\cup\{b\}\cup R_0(t')$ separates $z$ and $w$.
By Lemma  \ref{perterb-ray}, the rays $R_n(\alpha)$ and $R_n(\beta)$ land at the continuation $b_n$ of $b$ at $f_n$ for sufficiently large $n$. Then we get a sequence of simple curves $\G_n:=R_n(t)\cup\{b_n\}\cup R_n(t')$ for all large $n$.

In either case, according to Lemmas \ref{perterb-ray} and \ref{lem:internal-ray}, we have
$\limsup_{n\to\infty}\G_n=\G.$
Note that $\G$ separates $z,w$. Then $\G_n$ separates $z,w$ for all sufficiently large $n$. On the other hand, by taking a subsequence if necessary, we can assume that $R_{n}(\theta_n)$ is close to both $z$ and $w$ for large $n$.  It follows that there exist infinitely many $n$ for which $\G_n\cap R_{n}(\theta_n)\not=\emptyset$, a contradiction.
\end{proof}

\subsection{The limit of critical markings}\label{sec:not-closed}

\begin{proposition}\label{Julia-type-converge}
Let $f_n,n\geq 1,$ and $f$ be \pf polynomials in $\PPP_d$ such that $f_n\to f$ as $n\to\infty$, and let $\Theta_n$ be a critical marking of $f_n$  for each  large $n$.
If $(\Theta_n)$ Hausdorff converges to $\Theta$ as $n\to\infty$, then $\Theta$ is a weak critical marking of $f$.
 \end{proposition}

Before proving Proposition \ref{Julia-type-converge}, we show by two examples that the limit $\Theta$ is not necessarily a  critical marking of $f$. These examples are in fact the motivation for us to define the weak critical markings.

\begin{example}\label{E:not-conv}

We first consider the cubic polynomial $f(z)=z^3+3z^2/2$ (see Figure \ref{intersect}). Let $\{f_n,n\geq1\}$ be a sequence of postcritically-finite cubic polynomials in $\mathcal{P}_3$ converging to $f$, such that the rays with argument $1/3+\epsilon_n,2/3+\epsilon_n$ land at the unique Julia critical point of $f_n$, where $\epsilon_n\to 0^+$ as $n\to\infty$. We denote by $\G_n$ the union of the external rays of $f_n$ with arguments $1/3+\epsilon_n,2/3+\epsilon_n$, together with their common landing point. Thus, the deformation $U_n$ of $U$ is contained in the right-side component of $\C\setminus \G_n$, and hence $R_n(1/3)$ lands at $\partial U_n$ but $R_n(2/3)$ does not. As a consequence, we obtain a critical marking $\Theta_n$ of each $f_n$ as 
$$\Theta_n:=\big\{\mathcal{F}_n=\{0,1/3\},\mathcal{J}_n=\{1/3+\epsilon_n,2/3+\epsilon_n\}\big\}.$$
Clearly $\Theta_n$ Hausdorff converges to $\Theta=\big\{\{0,1/3\},\{1/3,2/3\}\big\}$, and $\Theta$ is a weak critical marking of $f$ but not a critical marking (since $R(1/3)$ is not left-supporting for $U$).

The second example is based on the cubic polynomial $f_{c_0}$ given in Figure \ref{portrait}, which admits two critical Julia markings $\big\{11/216, 83/216, 155/216\big\}$ and \big\{17/216, 89/216, 161/216\big\}.
Consider $\Theta=\big\{\Theta_1:=\{\frac{11}{216},\frac{83}{216}\}, \Theta_2:=\{\frac{89}{216},\frac{161}{216}\}\big\}$, a rational formal critical portrait of degree $3$ which is a weak critical marking, but not a critical marking, of $f_{c_0}$. The forward orbits of arguments in $\Theta$ are
\[\Theta_1\to \frac{11}{72}\to \frac{11}{24}\to\frac{3}{8}\rightleftarrows\frac{1}{8},\quad \Theta_2\to \frac{17}{72}\to \frac{17}{24}\to\frac{1}{8}\rightleftarrows\frac{3}{8}.\]
By perturbing $f_{c_0}$, one can find a sequence $(f_n)$ of \pf polynomials (with two distinct critical points) with $f_n \to f_{c_0}$, and such that each $f_n$ admits a critical marking of
the form
\[\Theta_n:=\big\{\Theta_{n,1}:=\Theta_1+\epsilon_{n,1},\Theta_{n,2}:=\Theta_2+\epsilon_{n,2} \big\},\]
with $\epsilon_{n,1},\epsilon_{n,2}\to 0$ as $n\to\infty$.
Then each $\Theta_n$ is a critical marking of $f_n$ and $\Theta_n \to \Theta$, but $\Theta$ is not a critical marking of $f_{c_0}$.
\end{example}

\begin{proof}[Proof of Proposition \ref{Julia-type-converge}]
Let $U$ be a critical Fatou component of $f$ and denote by $U_n$ the deformation of $U$ at $f_n$. 
Pick a critical marking for $f_n$, and let $\Theta(U_n)$ be the element associated to $U_n$ in this marking.

In the periodic case, each $\Theta(U_n)$ contains a unique periodic angle $\theta_n$ with period equal to that of $U_n$ and hence of $U$. By taking a subsequence if necessary, we can assume $\theta_n=\theta$ for large $n$.
Note that any $\Theta(U_n)$ is a subset of $\tau^{-1}(\tau(\theta))$ and $\#\Theta(U_n)={\rm deg}(f_n|_{U_n})={\rm deg}(f|_U)$ (by Lemma \ref{center-converge}), so we can further assume by taking a subsequence that $\Theta(U_n)$ is constant for large $n$, so we can write $\Theta(U_n) = \Theta$.
 According to Lemmas \ref{lem:internal-ray} and \ref{ray-convergence-2}, the rays with arguments in $\Theta$ land at the boundary of $U$. Furthermore, it follows from Lemma \ref{key} that the periodic angle in $\Theta$ is left-supporting for $U$.
In the strictly preperiodic case for $U$, by a similar argument and induction, we still get that $\Theta(U_n) = \Theta$ is constant  for large $n$, 
so that $\#\Theta ={\rm deg}(f|_U)$ and the external rays of $f$ with arguments in $\Theta$ land at the boundary of $U$.

Let $U_1,\ldots,U_s$ be all the critical Fatou components of $f$. The discussion above shows that the collection of sets $\{\Theta(U_1),\ldots,\Theta(U_s)\}$ is part of the critical marking $\Theta_n$ of $f_n$ for all large $n$ (by taking subsequences), and it is also a weak critical Fatou marking of $f$ as defined in Subsection \ref{sec:critical-Fatou}.

Now, we write each $\Theta_n$ as $\Theta_n:=\{\mathcal{F},\mathcal{L}_n\}$ with
$$\mathcal{F}:=\{\Theta(U_1),\ldots,\Theta(U_s)\} \text{\quad and\quad} \mathcal{L}_n:=\{\Theta_{n,1},\ldots,\Theta_{n,m}\},$$ for all large $n$, such that $\Theta_{n,i}\to \Theta_i$ as $n\to\infty$ and $1\leq i\leq m$. It follows immediately that $\#\Theta_{n,i}=\#\Theta_i$ for any $1\leq i\leq m$.
Note that each $\Theta_{n,i}$ corresponds to a critical point $c_{n,i}$ of $f_n$, and we can assume by taking subsequences that $c_{n,i}$ converge to a critical point $c_i$ of $f$, which must belong to $J_f$.

We claim that if all $c_{n,i}$ are in the Fatou set of $f_n$, the sequence of closed disks $(\ov{U(c_{n,i})})$ Hausdorff converges to the common landing point of the rays of $f_0$ with arguments in $\Theta_i$. It is enough to prove this result for any convergent subsequence of $(\ov{U(c_{n,i})})$, so we assume that $(\ov{U(c_{n,i})})$ converges in the Hausdorff metric to a connected compact set $S$. By Lemma \ref{ray-convergence-2}, the rays of $f_0$ with arguments in $\Theta_i$ land at $S$. So, to prove the claim, we only need to check that $S$ is a point.

By contradiction, we assume $\#S>1$ and choose $x\not=y\in S$. For all $n\geq1$, there exist $x_n\not=y_n\in \ov{U(c_{n,i})}$ such that $x_n\to x$ and $y_n\to y$ as $n\to\infty$.
As in the proof of Lemma \ref{ray-convergence-2},
we have $S\subset J_0$. Indeed, if there exists $z \in S$ which belongs to a Fatou component $U$, then by Lemma \ref{attracting-basin-converge}
it also belongs to its deformation $U_n$ for $n$ large, hence $U(c_{n,i})$ coincides with $U_n$, a contradiction. 

Let $[x,y]$ be the regulated arc in $K_f$ joining $x,y$.
 If $[x,y]\not\subset S$, then $[x,y]$ passes through a Fatou component $U$ of $f$. One can find two preperiodic rays $R_0(\alpha),R_0(\beta)$ of $f_0$ landing at $\partial U$ such that the orbits of their landing points avoid the critical points of $f_0$ and the set $R_0(\alpha)\cup U\cup R_0(\beta)$ separates $R_0(\theta)$ and $R_0(\eta)$. By Lemmas \ref{perterb-ray} and \ref{center-converge}, the rays $R_n(\alpha)$ and $R_n(\beta)$ of $f_n$ land at the boundary of $U_n$, the deformation of $U$ at $f_n$, and converge to $R_0(\alpha)$ and $R_0(\beta)$ respectively.  It implies that for all sufficiently large $n$, the simple curves $\G_n$, consisting of the union of $R_n(\alpha),R_n(\beta)$ and the internal rays in $U_n$ joining the landing points of $R_n(\alpha),R_n(\beta)$, separate $x_n$ and $y_n$.  Since $x_n$ and $y_n$ lie in the closure of the Fatou component $U(c_{n,i})$, this implies that $U_n = U(c_{n,i})$. This in turn yields that $c_i=\lim_{n\to\infty}c_{n,i}$ is a Fatou critical point, a contradiction to $c_i\in J_f$. 

If $[x,y]\subset S$ ($\subset J_f$), as what we shown in the proof of Lemma \ref{ray-convergence-2}, for each sufficiently large $n$, there exists a curve $\G_n:=R_n(\alpha)\cup\{z_n\}\cup R_n(\beta)$, consisting of two rays of $f_n$ and their common landing point, separating $x_n,y_n$, and these $\G_n$ converge to the simple curve $\G:=R_0(\alpha)\cup\{z\}\cup R_0(\beta)$ which separates $x$ and $y$.
 Therefore, for each sufficiently large $n$, we obtain a simple curve $\G_n$ which separates $x_n,y_n\in \ov{U(c_{n,i})}$ and is disjoint from $U(c_{n,i})$. This is impossible,
 so the claim is proven.

Now, to show $\Theta:=\{\FFF,\JJJ\}$ is a weak critical marking of $f$, we just need to check that $\JJJ$ satisfies properties (1)-(4) in the definition of weak critical Julia marking (see section \ref{sec:critical-julia-marking}). By the claim above and Lemma \ref{ray-convergence-2}, for each $i\in\{1,\ldots,m\}$, the rays of $f$ with arguments in $\Theta_i$ land at the common point $c_i$ of $f$. 
 On the other hand, the set $\{c_1,\ldots,c_m\}$ contains all the critical points of $f$ in the Julia set. To see this, note that any critical point $c\in J_f$  is an accumulation point of the critical points of $f_n$ according to Rouch\'{e}'s theorem. Furthermore, by Lemma \ref{center-converge}, the point $c$ cannot be an accumulation point
of the critical points of $f_n$ in the Fatou components corresponding to $\mathcal{F}$, hence it must be an accumulation point of $c_{n,1},\ldots,c_{n,m}, n\geq 1$. It follows that  $c\in\{c_1,\ldots,c_m\}$. The discussion above implies that properties (1)-(3) in the definition of weak critical Julia marking hold for $\JJJ$, so we just need to check property (4). Given a critical point $c\in J_f$, let $I_c$ denote the index set $$I_c:=\{i\in\{1,\ldots,m\}\mid c_i=c\}.$$
 Then $c_{n,i}\to c$ as $n\to\infty$ if and only if $i\in I_c$. Note that $c$ is a root of $f'(z)$ with multiplicity ${\rm deg}(f|_c)-1$. Then, by Rouch\'{e}'s Theorem for each sufficiently large $n$
 the function $f_n'$ has ${\rm deg}(f|_c)-1$ roots near $c$, counting with multiplicity. On the other hand, for each sufficiently large $n$, the points $c_{n,i}$ with $i\in I_c$ are exactly the roots of $f_n'$ near $c$, and each $c_{n,i}$ has multiplicity (\text{as a root of $f'_n$}) equal to
  $${\rm deg}(f_n|_{c_{n,i}})-1=\#\Theta(c_{n,i})-1=\#\Theta(c_i)-1.$$
It follows that the equation in property (4) holds.
 \end{proof}

 \subsection{The continuity of core entropy of polynomials}

 \begin{proof}[Proof of Theorem \ref{core-entropy-polynomial}]
 
Let $\Theta_n$ be a critical marking of the polynomial $f_n$. Since $f$ has only finitely many weak critical markings, by Proposition \ref{Julia-type-converge}, the sequence $(\Theta_n)$ can be subdivided into finitely many Hausdorff convergent subsequences. So it is enough to prove the theorem in the case that $(\Theta_n)$ Hausdorff converges to $\Theta$ as $n\to\infty$. Using Proposition \ref{Julia-type-converge}, the critical portrait $\Theta$ is a weak critical marking of $f$.
 Note that Theorem \ref{entropy-algorithm-1} shows that $h(\Theta_n)=h(f_n)$ and $h(\Theta)=h(f)$.
To complete the proof one needs to show $h(\Theta_n)\to h(\Theta)$ as $n\to\infty$. Let $m_n$ denote the primitive major induced by $\Theta_n$, and $m$ the primitive major induced by $\Theta$.  By Proposition \ref{pro:limit-of-major}.(2), the majors $m_n$ converge to $m$.
Moreover, by Lemma \ref{L:same-entro} and Lemma \ref{lem:entropy-continuity} one gets the equalities $h(\Theta_n) = h(m_n)$ and $h(\Theta) = h(m)$.
It follows from the claim above and Theorem \ref{theorem:main} that
$$h(\Theta_n)=h(m_n)\overset{{\rm Thm.}\ref{theorem:main}}{\longrightarrow} h(m)=h(\Theta),\  \text{as}\ n\to\infty.$$
The theorem is proven.

 \end{proof}

\end{document}